\documentclass[12pt]{amsart}

\usepackage{amssymb}
\usepackage{tabulary}
\usepackage[cmtip,all]{xy}

\newtheorem{theorem}{Theorem}[section]
\newtheorem{lemma}[theorem]{Lemma}
\newtheorem{proposition}[theorem]{Proposition}
\newtheorem{corollary}[theorem]{Corollary}

\theoremstyle{definition}
\newtheorem{definition}[theorem]{Definition}
\newtheorem{example}[theorem]{Example}

\theoremstyle{remark}
\newtheorem{remark}[theorem]{Remark}

\numberwithin{equation}{section}
\begin{document}

\title[Doubly warped structures]{Semi-Riemannian manifolds with a doubly warped structure}

\author[M. Guti\'errez]{Manuel Guti\'errez}
\address{Departamento de Algebra, Geometr\'ia y
Topolog\'ia\\
Universidad de M\'alaga. Spain.} \curraddr{}
\email{mgl@agt.cie.uma.es}
\thanks{This paper was supported in part by MEYC-FEDER Grant MTM2007-60016.}

\author[B. Olea]{Benjam\'in Olea}
\address{Departamento de Algebra, Geometr\'ia y
Topolog\'ia\\
Universidad de M\'alaga. Spain.} \curraddr{}
\email{benji@agt.cie.uma.es}

\subjclass[2000]{Primary 53C50, Secondary 53C12} \keywords{Quotient manifold, doubly warped product, doubly warped structure, decomposition
theorems.}
\date{}

\begin{abstract}
 We investigate manifolds obtained as a quotient of a doubly warped product. We show that they are always covered by the
 product of two suitable leaves. This allows us to prove, under regularity hypothesis, that these manifolds are a doubly warped product
 up to a zero measure subset formed by an union of leaves. We also obtain a necessary and sufficient condition which ensures the decomposition
 of the whole manifold and use it to give sufficient conditions of geometrical nature. Finally, we study the uniqueness of direct product
 decomposition in the nonsimply connected case.
\end{abstract}

\maketitle

\section{Introduction}

Let $(M_{i},g_{i})$ be two semi-Riemannian manifolds and
$\lambda_{i} :M_{1}\times M_{2}\rightarrow\mathbb{R}^{+}$ two
positive functions, $i=1,2$. The doubly twisted product
$M_{1}\times_{(\lambda_{1},\lambda_{2})}M_{2}$ is
the manifold $M_{1}\times M_{2}$ furnished with the metric $g=\lambda_{1}%
^{2}g_{1}+\lambda_{2}^{2}g_{2}$. If $\lambda_1\equiv 1$ or $\lambda_2\equiv 1$, then it is called a twisted product. On the other hand, when
$\lambda_{1}$ only depends on the second factor and $\lambda_{2}$ on the first one, it is called a doubly warped product (warped product if
$\lambda_1\equiv 1$ or $\lambda_2\equiv 1$). The most simple metric, the direct product, corresponds to the case $\lambda_1\equiv\lambda_2\equiv
1$.

Using the language of foliations, the classical De Rham-Wu Theorem says that two orthogonally, complementary and geodesic foliations (called a
direct product structure) in a complete and simply connected semi-Riemannian manifold give rise to a global decomposition as a direct product of
two leaves. This theorem can be generalized in two ways: one way is obtaining more general decompositions than direct products and the second
one is removing the simply connectedness hypothesis. The most general theorems in the first direction were obtained in \cite{Koike} and
\cite{Ponge}, where the authors showed that geometrical properties of the foliations determine the type of decomposition.

In the second direction, P. Wang obtained that a complete semi-Riemannian manifold furnished with a direct product structure is covered by the
direct product of two leaves, \cite{Wang}. Moreover, if the manifold is Riemannian, using the theory of bundle-like metrics, he showed that a
necessary and sufficient condition to obtain the global decomposition as a direct product is the existence of two regular leaves which intersect
each other at only one point. There are other remarkable works avoiding the simply connectedness but they treat with codimension one foliations,
see \cite{GutOle09} and references therein.

In this paper we study semi-Riemannian manifolds furnished with a \textit{doubly warped structure}, that is, two complementary, orthogonal and
umbilic foliations with closed mean curvature vector fields. If one of the mean curvature vector fields is identically null, then it is called a
\textit{warped structure}.

Doubly warped and warped structures appear in different situations. For example, Codazzi tensors with exactly two eigenvalues, one of them
constant \cite{Derdzinski}, Killing tensors \cite{Jelonek}, semi-Riemann submersions with umbilic fibres and some additional hypotheses
\cite{Escobales} and certain group actions \cite{Arouche,Lappas} lead to (doubly) warped structures. On the other hand, the translation of
physical content into geometrical language gives rise to warped structures which, by a topological simplification, are supposed global products.
In this way are constructed important spaces as Robertson-Walker, Schwarzschild, Kruskal, static spaces...

Manifolds with a doubly warped structure are locally a doubly warped product and, under completeness hypothesis, they are a quotient of a global
doubly warped product. This is why, after the preliminaries of section 2, we focus our attention on studying these quotients.

The main tool in this paper is Theorem \ref{teoprin}, which gives a normal semi-Riemannian covering map with a doubly warped product of two
leaves as domain. We use it to obtain a necessary and sufficient condition for a semi-Riemannian manifold with a doubly warped structure to be a
global product, which extends the one given in \cite{Wang} for direct products structures in the Riemannian setting. Other consequence of Theorem \ref{teoprin} is that any leaf is covered by a leaf without holonomy of the same foliation.

We study the space of leaves obtaining that, under regularity hypothesis, a manifold with a doubly warped structure is a fiber bundle over
the space of leaves. This allows us to compute the fundamental group of the space of leaves and to show that there is an open dense subset which
is isometric to a doubly warped product. We also give a result involving the curvature that ensures the global decomposition and apply it to
semi-Riemannian submersion with umbilic fibres.

As a consequence of the De Rham-Wu Theorem, it can be ensured the uniqueness of the direct product decomposition of a simply connected manifold
under nondegeneracy hypothesis. In the last section, we apply the decomposition results obtained to investigate the uniqueness of the
decomposition without the simply connectedness assumption.

\section{Preliminaries}

Given a product manifold $M_{1}\times M_{2}$ and $X\in\mathfrak{X}(M_{i})$, we will also denote $X$ to its elevation to
$\mathfrak{X}(M_{1}\times M_{2})$ and $P_{i}:TM_{1}\times TM_{2}\rightarrow TM_{i}$ will be the canonical projection. Unless it is explicitly
said, all manifolds are supposed to be semi-Riemannian. We write some formulaes about Levi-Civita connection and curvature, which are
established in \cite{Ponge}.

\begin{lemma}
\label{formulasderivadas}Let
$M_{1}\times_{(\lambda_{1},\lambda_{2})}M_{2}$ be a doubly twisted
product and call $\nabla^{i}$ the Levi-Civita connection of
$(M_{i},g_{i})$. Given $X,Y\in\mathfrak{X}(M_{1})$ and $V,W\in\mathfrak{X}%
(M_{2})$ it holds

\begin{enumerate}
\item
$\nabla_{X}Y=\nabla^{1}_{X}Y-g(X,Y)\nabla\ln\lambda_{1}+g(X,\nabla
\ln\lambda_{1})Y+g(Y,\nabla\ln\lambda_{1})X$.

\item
$\nabla_{V}W=\nabla^{2}_{V}W-g(V,W)\nabla\ln\lambda_{2}+g(V,\nabla
\ln\lambda_{2})W+g(W,\nabla\ln\lambda_{2})V$.

\item
$\nabla_{X}V=\nabla_{V}X=g(\nabla\ln\lambda_{1},V)X+g(\nabla\ln
\lambda_{2},X)V$.
\end{enumerate}
\end{lemma}

It follows that canonical foliations are umbilic and the mean curvature vector field of the first canonical foliation is
$N_{1}=P_{2}(-\nabla\ln\lambda_{1})$ whereas that of the second is $N_{2}=P_{1}(-\nabla\ln\lambda_{2})$.

\begin{lemma}
\label{curvatura} Let $M_{1}\times_{(\lambda_{1},\lambda_{2})}M_{2}$ be a doubly twisted product and take $\Pi_{1}=span(X,Y)$,
$\Pi_{2}=span(V,W)$ and $\Pi_{3}=span(X,V)$ nondegenerate planes where $X,Y\in TM_{1}$ and $V,W\in TM_{2}$ are unitary and orthogonal vectors.
Then the sectional curvature is given by

\begin{enumerate}
\item $K(\Pi_{1})=\frac{K^{1}(\Pi_{1})+g(\nabla\lambda_{1},\nabla\lambda_{1}%
)}{\lambda_{1}^{2}}-\frac{1}{\lambda_{1}} \Bigl(\varepsilon_{X}g(h_{1}%
(X),X)+\varepsilon_{Y}g(h_{1}(Y),Y)\Bigr).$

\item $K(\Pi_{2})=\frac{K^{2}(\Pi_{2})+g(\nabla\lambda_{2},\nabla\lambda_{2}%
)}{\lambda_{2}^{2}}-\frac{1}{\lambda_{2}} \Bigl(\varepsilon_{V}g(h_{2}%
(V),V)+\varepsilon_{W}g(h_{2}(W),W)\Bigr)$.

\item
$K(\Pi_{3})=-\frac{\varepsilon_{V}}{\lambda_{1}}g(h_{1}(V),V)-\frac
{\varepsilon_{X}}{\lambda_{2}}g(h_{2}(X),X)+\frac{g(\nabla\lambda_{1}
,\nabla\lambda_{2})}{\lambda_{1}\lambda_{2}}$,
\end{enumerate}
where $K^{i}$ is the sectional curvature of $(M_{i},g_{i})$,
$h_{i}$ is the hessian endomorphism of $\lambda_{i}$ and
$\varepsilon_{Z}$ is the sign of $g(Z,Z)$.
\end{lemma}

A vector field is called closed if its metrically equivalent one
form is closed.

\begin{lemma}
\label{lambdaenproducto}Let
$M_{1}\times_{(\lambda_{1},\lambda_{2})}M_{2}$ be a doubly twisted
product. It is a doubly warped product if and only if
$N_{i}=P_{3-i}(-\nabla\ln\lambda_{i})$ is closed, for $i=1,2$.
\end{lemma}
\begin{proof}
Suppose $d\omega_{1}=0$, where $\omega_{1}$ is the equivalent one form to $N_{1}$. If $X\in\mathfrak{X}(M_{1})$ and $V\in\mathfrak{X}(M_{2})$,
then $XV(\ln\lambda_{1})=-X\omega_{1}(V)=-d\omega_{1}(X,V)=0$. Thus there are functions $f_{1}\in C^{\infty}(M_{1})$ and $h_{1}\in
C^{\infty}(M_{2})$ such that $\lambda_{1}(x,y)=f_{1}(x)h_{1}(y)$ for all $(x,y)\in M_{1}\times M_{2}$. Analogously,
$\lambda_{2}(x,y)=f_{2}(x)h_{2}(y)$ for certain functions $f_{2}\in C^{\infty}(M_{1})$ and $h_{2}\in C^{\infty}(M_{2})$. Hence, taking conformal
metrics if necessary, $M_{1}\times_{(\lambda_{1},\lambda _{2})}M_{2}$ can be expressed as a doubly warped product. The \textit{only if} part is
trivial.
\end{proof}

We want to generalize the concept of doubly twisted or doubly
warped product to manifold which are not necessarily a topological
product.

\begin{definition}
\label{dts} Two complementary, orthogonal and umbilic foliations $(\mathcal{F}_{1},\mathcal{F}_{2})$ in a semi-Riemannian manifold is called a
doubly twisted structure. Moreover, if the mean curvature vectors of the foliations are closed, then it is called a doubly warped structure.
Finally, we say that it is a warped structure if one mean curvature vector is closed and the other one is zero.
\end{definition}

Notice that this last case is equivalent to one of the foliations being totally geodesic and the other one spherical, see \cite{Ponge} for the
definition.

We call $N_{i}$ the mean curvature vector field of $\mathcal{F}_{i}$ and  $\omega_{i}$, which we call \textit{mean curvature form}, to its
metrically equivalent one-form. The leaf of $\mathcal{F}_{i}$ through $x\in M$ is denoted by $F_{i}(x)$ and $\mathcal{F}_{i}(x)$ will be the
tangent plane of $F_{i}(x)$ at the point $x$. If there is not confusion or if the point is not relevant, we simply write $F_{i}$. We always put
the induced metric on the leaves.

\begin{remark}\label{remarkglobal}If $M$ has a doubly twisted (warped) structure, then we can take around any point an adapted chart to both foliations. Lemma
\ref{lambdaenproducto} and Proposition 3 of \cite{Ponge} show that M is locally isometric to a doubly twisted (warped) product. In the doubly
warped structure case, the condition on the mean curvature vectors in Theorem 5.4 of \cite{Koike} can be easily checked. So, if the leaves of
one of the foliations are complete we can apply this theorem to obtain that $M$ is a quotient of a global doubly warped product.
\end{remark}

Given a curve $\alpha:[0,1] \rightarrow M$ we call
$\alpha_{t}:[0,t]\rightarrow M$, $0\leq t\leq1$, its restriction.

\begin{definition}
\label{Traslacion Adaptada}Let $M$ be a semi-Riemannian manifold
with $\mathcal{F}_{1}$ and $\mathcal{F}_{2}$ two orthogonal and
complementary foliations. Take $x\in M$, $v\in\mathcal{F}_{2}(x)$
and $\alpha :[0,1] \rightarrow F_{1}(x)$ a curve with
$\alpha(0)=x$. We define the adapted translation of $v$ along
$\alpha_{t}$ as $A_{\alpha_{t}}(v)=\exp\left(
-\int_{\alpha_{t}}\omega_{2}\right)  W(t)$, where $W$ is the normal
parallel translation to $\mathcal{F}_1$ of $v$ along
$\alpha$, \cite{ONeill}.
\end{definition}

In the same way we can define the adapted translation of a vector
of $\mathcal{F}_{1}(x)$ along a curve in $F_{2}(x)$. Observe that
$|A_{\alpha_{t}}(v)|=|v|\exp\left(
-\int_{\alpha_{t}}\omega_{2}\right)  $.

\begin{lemma}
\label{translationinldws} Let $M=M_{1}\times M_{2}$ be a
semi-Riemannian manifold such that the canonical foliations
constitute a doubly twisted structure. Take $\alpha:[0,1]
\rightarrow M_{1}$ a curve with $\alpha(0)=a$ and $v_{b}\in
T_{b}M_{2}$. The adapted translation of $(0_{a},v_{b})$ along the
curve $\gamma(t)=(\alpha(t),b)$ is
$A_{\gamma_{t}}(0_{a},v_{b})=(0_{\alpha (t)},v_{b})$.
\end{lemma}
\begin{proof} First we show that formula 3 of Lemma \ref{formulasderivadas}
is still true in this case. In fact, take $X,Y\in\mathfrak{X}
(M_{1})$ and $V,W\in\mathfrak{X}(M_{2})$. Since $[X,V]=0$ we have
\begin{equation*}
g(\nabla_{X}V,W)=-g(X,\nabla_{V}W)=-g(V,W)g(X,N_{2}),
\end{equation*}
and analogously $g(\nabla_{X}V,Y)=-g(V,N_{1})g(X,Y)$. Therefore,
it follows that $\nabla_{X}V=-\omega_{1}(V)X-\omega_{2}(X)V$ for
all $X\in\mathfrak{X}(M_{1})$ and $V\in\mathfrak{X}(M_{2})$.

Take $V(t)=(0_{\alpha(t)},v_{b})$ and $W(t)=\lambda(t)V(t)$, where
$\lambda(t)=\exp\left(  \int_{\gamma_{t}}\omega_{2}\right)$. We
only have to check that $W(t)$ is the normal parallel translation
of $(0_{a},v_{b})$ along $\gamma$. But this is an immediate
consequence of
\begin{align*}
\frac{DW}{dt} &  =\lambda^{\prime}V+\lambda(-\omega_{1}%
(V)\gamma^{\prime}-\omega_{2}(\gamma^{\prime})V) =-\lambda\omega_{1}(V)\gamma^{\prime}.
\end{align*}
\end{proof}

Given a foliation $\mathcal{F}$ we call $Hol(F)$ the holonomy group of a leaf $F$ (see \cite{Camacho} for definitions and properties). We say
that $F$ has not holonomy if its holonomy group is trivial. The foliation $\mathcal{F}$ has not holonomy if any leaf has not holonomy.

\begin{lemma}
\label{TraslAdaptada-AplicHolon}Let $M$ be a semi-Riemannian
manifold with $(\mathcal{F}_{1},\mathcal{F}_{2})$ a doubly twisted
structure. Take $x\in M$ and $\alpha:[0,1] \rightarrow F_{1}(x)$ a
loop at $x$. If $f\in Hol(F_{1}(x))$ is the holonomy map
associated to $\alpha$, then $f_{\ast_{x}}(v)=A_{\alpha }(v)$ for
all $v\in\mathcal{F}_{2}(x).$
\end{lemma}
\begin{proof} It is sufficient to show it locally. Take an open set of $x$
isometric to a doubly twisted product
$U_{1}\times_{(\lambda_{1},\lambda_{2})}U_{2}$ where $U_{i}$ is an
open set of $F_{i}(x)$ with $x\in U_{i}$. If $\alpha(t)\in U_{1} $
, for $0\leq t\leq t_{0},$ then the holonomy map associated to
this arc is $f:\{x\}\times U_{2}\rightarrow\{\alpha(t_{0})\}\times
U_{2}$ given by
$f(x,y)=(\alpha(t_{0}),y)$ and clearly $f_{\ast_{x}}(0_{x},v_{x}%
)=(0_{\alpha(t_{0})},v_{x})=A_{\alpha_{t_{0}}}(v_{x})$.
\end{proof}

Observe that, in the doubly warped structure case, holonomy maps are homotheties and thus
 they are determined by their derivative at a point. Therefore, a leaf $F_1$ has not holonomy
if and only if for any loop $\alpha$ in $F_1$ it holds $A_\alpha=id$.

We finish this section relating two doubly twisted structure via a
local isometry.

\begin{lemma}\label{dosdts} Let $\overline{M}$ and $M$ be two semi-Riemannian
manifold with
a doubly twisted structure $(\overline{\mathcal{F}}_{1},\overline{\mathcal{F}%
}_{2})$ and $(\mathcal{F}_{1},\mathcal{F}_{2})$ respectively. Take
$f:\overline{M}\rightarrow M$ a local isometry which preserves the
doubly
twisted structure, that is, $f_{\ast_{x}}(\overline{\mathcal{F}}%
_{i}(x))=\mathcal{F}_{i}(f(x))$ for all $x\in\overline{M}$ and
$i=1,2$. Then $f$ also preserves
\begin{enumerate}
\item the leaves, $f(\overline{F}_{i}(x))\subset F_{i}(f(x))$.
\item the mean curvature vector fields, $f_{\ast_{x}}(\overline{N}%
_{i}(x))=N_{i}(f(x))$. \item the mean curvature forms,
$f^{\ast}(\omega_{i})=\overline{\omega}_{i}$. \item the adapted
translation, that is, if $\gamma:[0,1] \rightarrow\overline
{F}_{1}(x)$ is a curve with $\gamma(0)=x$ and $v\in\overline{\mathcal{F}}%
_{2}(x)$, then $A_{f\circ\gamma_{t}}\left(  f_{\ast_{x}}(v)\right)
=f_{\ast_{\gamma(t)}}\left(  A_{\gamma_{t}}(v)\right)$ and
analogously for a curve in $\overline{F}_{2}(x)$.
\end{enumerate}
\end{lemma}
\begin{proof}
\begin{enumerate}
\item It is immediate.
\item Take $\overline{X},\overline{Y}\in
\overline{\mathcal{F}}_{1}$ and call $X=f_{\ast}(\overline{X})$
and $Y=f_{\ast}(\overline{Y})$ in a suitable open set. Since $f$
is a local isometry, $g(X,Y)N_1=P_2( \nabla_{X}Y) =
f_{\ast}(\overline{P}_2( \overline{\nabla}_{\overline{X}}\overline
{Y}))=g(X,Y)f_{\ast}(\overline{N}_{1})$. Hence
$N_{1}=f_{\ast}({\overline{N}_{1}})$ and analogously for
$\overline{N}_{2}$.

\item Immediate from point (2).

\item Take $\overline{W}(t)$ the parallel translation of $v$ normal to $\overline{\mathcal{F}}_1$ along $\gamma$. Since $f$ is a local isometry
which preserves the foliations, $f_{\ast}(\overline{W}(t))$ is the normal parallel translation of $f_{\ast }(v)$ along $f\circ\gamma$. Using
point (3),
 $\int_{\gamma
_{t}}\overline{\omega}_{2}=\int_{f\circ\gamma_{t}}\omega_{2}$ and
therefore $A_{f\circ\gamma_{t}}\left( f_{\ast_{x}}(v)\right)
=f_{\ast _{\gamma(t)}}\left( A_{\gamma_{t}}(v)\right)$.
\end{enumerate}
\end{proof}

\section{Quotient of a doubly warped product}
From now on, $M_{1}\times_{(\lambda_{1},\lambda_{2})}M_{2}$ will be a doubly warped product and $\Gamma$ a group of isometries such that:

\begin{enumerate}
\item $\Gamma$ acts in a properly discontinuous manner. \item It preserves the canonical foliations. This implies that if $f\in\Gamma $, then
$f=\phi\times\psi:M_{1}\times M_{2}\rightarrow M_{1}\times M_{2}$, where $\phi:M_{1}\rightarrow M_{1}$ and $\psi:M_{2}\rightarrow M_{2}$ are
homotheties with factor $c_{1}^2$ and $c_{2}^2$ respectively, such that $\lambda_{1}\circ
\psi=\frac{1}{c_{1}}\lambda_{1}$ and $\lambda_{2}\circ\phi=\frac{1}{c_{2}%
}\lambda_{2}$.
\end{enumerate}

The semi-Riemannian manifold
$M=\left(M_{1}\times_{(\lambda_{1},\lambda _{2})}M_{2}\right)
/\Gamma$ has a doubly warped structure, which, as always, we call
$(\mathcal{F}_{1},\mathcal{F}_{2})$. We are going to work with
$\mathcal{F}_{1}$ because all definitions and results are stated
analogously for $\mathcal{F}_{2}$.

If we take the canonical projection $p:M_{1}\times
M_{2}\rightarrow M$, which is a semi-Riemannian covering map,
applying Lemma \ref{dosdts} we have $p(M_{1}\times\{b\})\subset
F_{1}(p(a,b))$ for all $(a,b)\in M_{1}\times M_{2}$. We call
$p_{1}^{(a,b)}:M_{1}\times\{b\}\rightarrow F_{1}(p(a,b))$ the
restriction of $p$.

\begin{lemma}\label{coveringinducidos} Let $M=\left(M_{1}\times_{(\lambda_{1},\lambda _{2})}M_{2}\right)  /\Gamma$ be
a quotient of a doubly warped product. We take $p:M_1\times
M_2\rightarrow M$ the canonical projection, $(a,b)\in M_{1}\times
M_{2}$ and $x=p(a,b)$. Then, the restriction
$p_{1}^{(a,b)}:M_{1}\times\{b\}\rightarrow F_{1}(x)$ is a normal
semi-Riemannian covering map.
\end{lemma}
\begin{proof} It is clear that $p_{1}^{(a,b)}:M_{1}\times\{b\}\rightarrow
F_{1}(x)$ is a local isometry. Let $\gamma:[0,1] \rightarrow
F_{1}(x)$ be a curve with $\gamma(0)=x$. Since $p:M_{1}\times
M_{2}\rightarrow M$ is a covering map, there is a lift
$\alpha:[0,1] \rightarrow M_{1}\times M_{2}$ with $\alpha
(0)=(a,b)$. But
$p_{\ast}(\alpha^{\prime}(t))=\gamma^{\prime}(t)\in
\mathcal{F}_{1}(\gamma(t))$ and $p$ preserves the foliations, so
$\alpha(t)$ is a curve in $M_{1}\times\{b\}$. Applying Theorem 28
of \cite[pg. 201]{ONeill}, we get that $p_{1}^{(a,b)}$ is a
covering map. Now we show that it is normal. Take $a^{\prime }\in
M_{1}$ such that $p_{1}^{(a,b)}(a^{\prime},b)=x$. Then, there
exists $f\in\Gamma$ with $f(a^{\prime},b)=(a,b)$ and since $f$
preserves the canonical foliations,
$f(M_{1}\times\{b\})=M_{1}\times\{b\}$. So, the restriction of $f$
to $M_{1}\times\{b\}$ is a deck transformation of the covering
$p_{1}^{(a,b)}$ which sends $(a^{\prime},b)$ to $(a,b)$.
\end{proof}

Let $\Gamma_{1}^{(a,b)}$ be the group of deck transformations of $p_{1}^{(a,b)}:M_{1}\times\{b\}\rightarrow F_{1}(p(a,b))$. When there is not
confusion with the chosen point, we simply write $p_{1}$ and $\Gamma_{1}$ instead of $p_{1}^{(a,b)}$ and $\Gamma_{1}^{(a,b)}$.

If $\phi\in\Gamma_{1}$, in general, it does not exist
$\psi\in\Gamma_{2}$ with $\phi\times\psi\in\Gamma$ and it has not
to hold that $\lambda_{2}\circ\phi=\frac{1}{c_{2}}\lambda_{2}$, as
it was said at the beginning of this section.

\begin{lemma}\label{subgrupo} Let $M=\left(  M_{1}\times_{(\lambda_{1},\lambda_{2})}%
M_{2}\right)  /\Gamma$ be a quotient of a doubly warped product.
Fix $(a_{0},b_{0})\in M_{1}\times M_{2}$ and
$x_{0}=p(a_{0},b_{0})$ such that the leaf $F_{1}(x_{0})$ has not
holonomy. Then
\begin{enumerate}
\item $\lambda_{2}\circ\phi=\lambda_{2}$ for all $\phi\in\Gamma_{1}%
^{(a_{0},b_{0})}$. \item $\phi\times id\in\Gamma$ for all
$\phi\in\Gamma_{1}^{(a_{0},b_{0})}$ and$\
\Gamma_{1}^{(a_{0},b_{0})}\times\{id\}$ is a normal subgroup of
$\Gamma$.
\end{enumerate}
\end{lemma}
\begin{proof}
\begin{enumerate}
\item The mean curvature forms $\overline{\omega}_{2}$ and $\omega_{2}$ of the foliations in $M_{1}\times M_{2}$ and $M$ respectively are
closed. Thus every point in $M$ has an open neighborhood where $\omega _{2}=df_{2}$ for certain function $f_2$, and analogously every point in
$M_{1}\times M_{2}$ has an open neighborhood where $\overline{\omega}_{2}=d\overline{f}_{2}$, being in this case
$\overline{f}_{2}=-\ln\lambda_{2}$. Since $p_{1}\circ\phi=p_{1}$, we have $\phi^{\ast}(p^{\ast}(\omega_{2}))=p^{\ast}(\omega_{2})$ and therefore
$f_{2}\circ p\circ\phi=f_{2}\circ p+k_{1}$ for certain constant $k_{1}$. On the other hand, using Lemma \ref{dosdts}, we have $f_{2}\circ p=-\ln
\lambda_{2}+k_{2}$ for some constant $k_{2}$.

Joining the last two equations, we get
$\lambda_{2}\circ\phi=c\lambda_{2}$ for certain constant $c$. This
formula must be true in the whole $M_{1}$ with the same constant
$c$ and, in fact,
$c=\frac{\lambda_{2}(a_{1})}{\lambda_{2}(a_{0})}$ where
$a_{1}=\phi(a_{0})$. Take $\alpha:[0,1] \rightarrow M_{1}
\times\{b_{0}\}$ with $\alpha(0)=(a_{0},b_{0})$,
$\alpha(1)=(a_{1},b_{0})$ and $w\in T_{b_{0}}M_{2}$ a non
lightlike vector. Using Lemma \ref{translationinldws} we have
$A_{\alpha}(0_{a_{0}},w)=(0_{a_{1}},w)$ and since $p$ preserves
the adapted translation,
\begin{equation*}
A_{p\circ\alpha}\Bigl(p_{\ast_{(a_{0},b_{0})}}(0_{a_{0}},w)\Bigr)=p_{\ast
_{(a_{1},b_{0})}}(0_{a_{1}},w).
\end{equation*}

Using that $F_{1}(x_{0})$ has not holonomy and Lemma \ref{TraslAdaptada-AplicHolon}, we obtain
\begin{equation*}
p_{\ast_{(a_{1},b_{0})}}(0_{a_{1} },w)=p_{\ast_{(a_{0},b_{0})}}(0_{a_{0}},w),
\end{equation*}
taking norms, we get $c=1$.

\item Take $\phi\in\Gamma_{1}$. Since
$\lambda_{2}\circ\phi=\lambda_{2}$, it follows that $\phi\times
id$ is an isometry of
$M_{1}\times_{(\lambda_{1},\lambda_{2})}M_{2}$. Now, to show that
$p\circ\left(  \phi\times id\right) =p$ it is
enough to check $(p\circ\phi\times id)_{\ast(a_{0},b_{0})}=p_{\ast(a_{0}%
,b_{0})}$. Using $p_{1}\circ\phi=p_{1}$, we see that
\begin{equation*}
(p\circ\left(  \phi\times id\right)  )_{\ast(a_{0},b_{0})}(v,0_{b_{0}%
})=p_{\ast(a_{0},b_{0})}(v,0_{b_{0}})
\end{equation*}
for all $v\in T_{a_{0}}M_{1}$. Given $w\in T_{b_{0}}M_{2}$, we
take $\alpha:[0,1] \rightarrow M_{1}\times\{b_{0}\}$ a curve from
$(a_{0},b_{0})$ to
$(a_{1},b_{0})$, where $\phi(a_{0})=a_{1}$. Then $A_{\alpha}(0_{a_{0}%
},w)=(0_{a_{1}},w)$ and
\begin{align*}
(p\circ(\phi\times id))_{\ast_{(a_{0},b_{0})}}(0_{a_{0}},w)
=p_{\ast _{(a_{1},b_{0})}}(0_{a_{1}},w)=\\
p_{\ast_{(a_{1},b_{0})}}\Bigl(A_{\alpha}(0_{a_{0}},w)\Bigr)
=A_{p\circ
\alpha}\Bigl(p_{\ast_{(a_{0},b_{0})}}(0_{a_{0}},w)\Bigr)=p_{\ast_{(a_{0},b_{0})}}(0_{a_{0}},w),
\end{align*}
where the last equality holds because $F_{1}(x_{0})$ has not
holonomy. Therefore $\phi\times id\in\Gamma$ and
$\Gamma_{1}\times\{id\}$ is a subgroup of $\Gamma$.

Now to prove that it is a normal subgroup, we take
$\phi\in\Gamma_{1}$ and $f\in\Gamma$ and show that
$f^{-1}\circ\left(  \phi\times id\right)  \circ f\in\Gamma_{1}$.
Since $f$ preserves the foliations, $f^{-1}\circ\left(
\phi\times id\right)  \circ f$ takes $M_{1}\times\{b_{0}\}$ into $M_{1}%
\times\{b_{0}\}$ and therefore we can consider $h=\left.
f^{-1}\circ\left(
\phi\times id\right)  \circ f\right\vert _{M_{1}\times\{b_{0}\}}\in\Gamma_{1}%
$. But $h\times id$ coincides with $f^{-1}\circ\left(  \phi\times
id\right) \circ f$ at $a_{0}$, thus $f^{-1}\circ\left(  \phi\times
id\right)  \circ f=h\times id\in\Gamma_{1}\times\{id\}$.
\end{enumerate}
\end{proof}

\begin{remark}\label{subgrupo2}Observe that we have used that $F_1(x_0)$ has not holonomy only to ensure $A_{p\circ
\alpha}=id$, thus we have a little bit more general result. Suppose that $\gamma:[0,1]\rightarrow F_1(x)$ is a loop at $x\in M$ such that its
associated holonomy map is trivial. Take $\alpha:[0,1]\rightarrow M_1\times \{b\}$ a lift through $p_1:M_1\times \{b\}\rightarrow F_1(x)$ with
basepoint $(a,b)$ and suppose $\alpha(1)=(a',b)$. If $\phi\in\Gamma_1$ with $\phi(a)=a'$, then it can be proven, identically as in the above
lemma, that for this deck transformation it holds $\lambda_2\circ\phi=\lambda_2$ and $\phi\times id\in\Gamma$. This will be used in Theorem
\ref{Covering hojas}.
\end{remark}

Now, we can give the following theorem, which is the main tool of
this paper (compare with Theorem 2 of \cite{Wang} and Theorem 7 of
\cite{ShaZhuIgo}).

\begin{theorem}\label{teoprin} Let $M=\left(  M_{1}\times_{(\lambda_{1},\lambda_{2})}%
M_{2}\right)  /\Gamma$ be a quotient of a doubly warped product.
Fix $(a_{0},b_{0})\in M_{1}\times M_{2}$ and
$x_{0}=p(a_{0},b_{0})$. The leaf $F_{1}(x_{0})$ has not holonomy
if and only if there exists a semi-Riemannian
normal covering map $\Phi:F_{1}(x_{0})\times_{(\lambda_{1},\rho_{2})}%
M_{2}\rightarrow M$, where
$\rho_{2}:F_{1}(x_{0})\rightarrow\mathbb{R}^{+}$. Moreover, the
following diagram is commutative
\begin{equation*}
\xymatrix{ M_1\times M_2\ar[d]_{p_1^{(a_0,b_0)}\times id} \ar[r]^p
& M \\ F_1(x_0)\times M_2 \ar[ru]_{\Phi} }
\end{equation*}
In particular, $\Phi(x,b_{0})=x$ for all $x\in F_{1}(x_{0})$.
\end{theorem}
\begin{proof} Suppose that $F_{1}(x_{0})$ has not holonomy. Since $\Gamma_{1}\times \{id\}$
is a normal subgroup of $\Gamma$, there exists a normal covering
map
\begin{equation*}
\Phi:\left(  M_{1}\times_{(\lambda_{1},\lambda_{2})}M_{2}\right) /\left( \Gamma_{1}\times \{id\}\right)  \rightarrow M.
\end{equation*}
But $\left(  M_{1}\times_{(\lambda_{1},\lambda_{2})}M_{2}\right) /\left( \Gamma_{1}\times \{id\}\right)  $ is isometric to $F_{1}(x_{0})\times
_{(\lambda_{1},\rho_{2})}M_{2}(x_{0})$ for certain function $\rho_{2}$ with $\rho_{2}\circ p_{1}=\lambda_{2}$ and by construction
$\Phi\circ\left( p_{1}\times id\right)  =p.$

Conversely, we suppose the existence of such semi-Riemannian
covering. Take $\alpha:[0,1] \rightarrow F_{1}(x_{0})$ a loop at
$x_{0}$, $w\in
T_{x_{0}}F_{2}(x_{0})$ and $v\in T_{b_{0}}M_{2}$ with $p_{\ast_{(a_{0},b_{0}%
)}}(0,v)=w$. Then
\begin{equation*}
\Phi_{\ast_{(x_{0},b_{0})}}(0,v)=\Phi_{\ast_{(x_{0},b_{0})}}\Bigl(\left(
p_{1}\times id\right)  _{\ast_{(a_{0},b_{0})}}(0,v)\Bigr)=p_{\ast
_{(a_{0},b_{0})}}(0,v)=w.
\end{equation*}

Now, using Lemma \ref{translationinldws}, \ref{dosdts} and that the holonomy in $F_1(x_0)\times_{(\lambda_1,\rho_2)} M_2$ is trivial,
\begin{equation*}
A_{\alpha}(w)=\Phi_{\ast_{(x_{0},b_{0})}}\left(  A_{(\alpha,b_{0}%
)}(0,v)\right)  =\Phi_{\ast_{(x_{0},b_{0})}}(0,v)=w
\end{equation*}
and therefore $F_{1}(x_{0})$ has not holonomy.
\end{proof}

We say that $x_{0}\in M$ has not holonomy if $F_{1}(x_{0})$ and $F_{2}(x_{0}) $ have not holonomy.

\begin{corollary}\label{corprin} Let $M=\left(  M_{1}\times_{(\lambda_{1},\lambda_{2})}%
M_{2}\right)  /\Gamma$ be a quotient of a doubly warped product.
Fix $(a_{0},b_{0})\in M_{1}\times M_{2}$ and
$x_{0}=p(a_{0},b_{0})$. The point $x_{0}$ has not holonomy if and
only if there is a semi-Riemannian normal
covering map
\begin{equation*}
\Phi:F_{1}(x_{0})\times_{(\rho_{1},\rho_{2})}F_{2}%
(x_{0})\rightarrow M,
\end{equation*}
where $\rho_{1}:F_{2}(x_{0})\rightarrow\mathbb{R}^{+}$ and $\rho_{2}:F_{1}(x_{0})\rightarrow\mathbb{R}^{+}$. Moreover, the following diagram is
commutative
\begin{equation*}
\xymatrix{ M_1\times M_2 \ar[d]_{p_1^{(a_0,b_0)}\times
p_2^{(a_0,b_0)}} \ar[r]^{p} & M \\ F_1(x_0)\times F_2(x_0)
\ar[ur]_{\Phi} }
\end{equation*}
In particular, $\Phi(x,x_{0})=x$ and $\Phi(x_{0},y)=y$ for all
$x\in F_{1}(x_{0})$ and $y\in F_{2}(x_{0})$.
\end{corollary}

\begin{remark} It is known that the set of leaves without holonomy is dense on
$M$, \cite{Tischler}. Thus, we can always take a point $x_{0}\in M$ without holonomy and apply Corollary \ref{corprin}.
\end{remark}

\begin{theorem} Let $M=\left(  M_{1}\times_{(\lambda_{1},\lambda_{2})}M_{2}\right)
/\Gamma$ be a quotient of a doubly warped product. If $x_{0}\in M$
has not holonomy then
\begin{equation*}
card(F_{1}(x_{0})\cap F_{2}(x_{0}))=card(\Phi^{-1}(x_{0})),
\end{equation*}
where $\Phi:F_{1}(x_{0})\times_{(\rho_{1},\rho_{2})}F_{2}(x_{0})\rightarrow M $ is the semi-Riemannian covering map of Corollary \ref{corprin}.
\end{theorem}

\begin{proof} Recall that, by construction, $\Phi(x,y)=p(a,b)$ where $a\in
M_{1}$ with $p(a,b_{0})=x$ and $b\in M_{2}$ with $p(a_{0},b)=y$.
Now we define $\Lambda:\Phi^{-1}(x_{0})\rightarrow
F_{1}(x_{0})\cap F_{2}(x_{0})$ by $\Lambda(x,y)=x$. First we show
that $\Lambda$ is well defined. If $(x,y)\in\Phi^{-1}(x_{0})$ then
$p(a,b)=x_{0}$, where $p(a,b_{0})=x$ and
$p(a_{0},b)=y$. Hence $p(\{a\}\times M_{2})\subset F_{2}(p(a,b))=F_{2}%
(x_{0})$ and $p(M_{1}\times\{b_{0}\})\subset F_{1}(p(a_{0},b_{0}))=F_{1}%
(x_{0})$, thus
\begin{equation*}
x=p(a,b_{0})=p(M_{1}\times\{b_{0}\}\cap\{a\}\times M_{2})\in F_{1}(x_{0})\cap F_{2}(x_{0}).
\end{equation*}

Now we check that $\Lambda$ is onto. Take $x\in F_{1}(x_{0})\cap
F_{2}(x_{0}) $. Since $p_{1}:M_{1}\times\{b_{0}\}\rightarrow
F_{1}(x_{0})$ is a covering map there exists $a\in M_{1}$ such
that $p(a,b_{0})=x$. But $p_{2}^{(a,b_{0})}:\{a\}\times
M_{2}\rightarrow F_{2}(x)=F_{2}(x_{0})$ is a covering map too,
therefore there is $(a,b)\in \{a\}\times M_{2}$ such that
$p(a,b)=x_{0}$. If we call $y=p(a_{0},b)$, then $\Phi(x,y)=x_{0}$
and $\Lambda(x,y)=x$.

Finally, we show that $\Lambda$ is injective. Take $(x,y),(x,y^{\prime}%
)\in\Phi^{-1}(x_{0})$ and $a\in M_{1}$, $b,b^{\prime}\in M_{2}$
such that $p(a,b_{0})=x$, $p(a_{0},b)=y$ and
$p(a_{0},b^{\prime})=y^{\prime }$. Consider the covering
$p_{2}^{(a,b_{0})}:\{a\}\times M_{2}\rightarrow F_{2}(x_{0})$.
Since $p_{2}^{(a,b_{0})}(a,b)=p_{2}^{(a,b_{0})}(a,b^{\prime})$ and
this covering is normal, there exist a deck transformation
$\psi\in\Gamma_{2}^{(a,b_{0})}$ such that
$\psi(a,b)=(a,b^{\prime})$. But $F_{2}(x_{0})$ has not holonomy,
so Lemma \ref{subgrupo} assures that $id\times\psi\in\Gamma$. Now,
$(id\times \psi)(a_{0},b)=(a_{0},b^{\prime})$ and thus
$y=y^{\prime}$.
\end{proof}

Now we give a necessary and sufficient condition for a doubly
warped structure to be a global doubly warped product, which
extends the one given in \cite{Wang} for direct products and
Riemannian manifolds.

\begin{corollary}\label{descglobal} Let $M=\left(  M_{1}\times_{(\lambda_{1},\lambda_{2})}%
M_{2}\right)  /\Gamma$ be a quotient of a doubly warped product
and $x_{0}\in M$. Then $M$ is isometric to the doubly warped product $F_{1}(x_{0}%
)\times_{(\rho_{1},\rho_{2})}F_{2}(x_{0})$ if and only if $x_{0}$
has not holonomy and $F_{1}(x_{0})\cap F_{2}(x_{0})=\{x_{0}\}$.
\end{corollary}

Condition $F_{1}(x_{0})\cap F_{2}(x_{0})=\{x_{0}\}$ alone is not
sufficient to split $M$ as a product $F_{1}(x_{0})\times
F_{2}(x_{0})$, as intuition perhaps suggests. The M\"{o}bius trip
illustrates this point.

\begin{theorem} Let $M=\left(  M_{1}\times_{(\lambda_{1},\lambda_{2})}M_{2}\right)
/\Gamma$ be a quotient of a doubly warped product. If $x_{0}$ has
not holonomy then
\begin{equation*}
card(F_{1}(x)\cap F_{2}(x))\leq card(F_{1}(x_{0})\cap
F_{2}(x_{0}))
\end{equation*}
for all $x\in M$.
\end{theorem}
\begin{proof} Take $(a_{0},b_{0})\in M_{1}\times M_{2}$ such that
$p(a_{0},b_{0})=x_{0}$. First suppose that $x\in F_{1}(x_{0})$ and $F_{1}(x)\cap F_{2}(x)=\{x_{i}:i\in I\}$. If we take $a\in M_{1}$ such that
$p(a,b_{0})=x$, then we know that $p_{2}^{(a,b_{0})}:\{a\}\times M_{2}\rightarrow F_{2}(x)$ is a covering map, so we can take $b_{i}\in M_{2}$
with $p_{2}^{(a,b_{0})}(a,b_{i})=x_{i}$. If we call $y_{i}=p_{2}^{(a_{0},b_{0})}(a_{0},b_{i})$, then both $x_{i},y_{i}\in
p(M_{1}\times\{b_{i}\})=F_{1}(p(a,b_{i}))=F_{1}(x_{0})$, and moreover, since $y_{i}\in F_{2}(x_{0})$ we have $y_{i}\in F _{1}(x_{0})\cap
F_{2}(x_{0})$.

Now, we show that the map $\Lambda:F_{1}(x)\cap
F_{2}(x)\rightarrow F_{1}(x_{0})\cap F_{2}(x_{0})$ given by
$\Lambda(x_{i})=y_{i}$ is injective. If
$y_i=p_{2}^{(a_{0},b_{0})}(a_{0},b_{i})=p_{2}^{(a_{0},b_{0})}(a_{0},b_{j})=y_j$
for
$i\neq j$ then there is $\psi\in\Gamma_{2}^{(a_0,b_0)}$ such that $\psi(a_{0}%
,b_{i})=(a_{0},b_{j})$. Since $F_{2}(x_{0})$ has not holonomy
(Lemma \ref{subgrupo}), $id\times\psi\in\Gamma$ and it sends
$(a,b_{i})$ to $(a,b_{j})$. Therefore $x_{i}=x_{j}$. This shows
that $card(F_{1}(x)\cap F_{2}(x))\leq card(F_{1}(x_{0})\cap
F_{2}(x_{0}))$ when $x\in F_1(x_0)$.

Take now an arbitrary point $x\in M$ and $(a,b)\in M_{1}\times
M_{2}$ with $p(a,b)=x$. We have that $F_{2}(x)$ intersects $F
_{1}(x_{0})$ at some point $z=p(a,b_{0})$. In the same way as
above, using that $F_{1}(x_{0})$ has not holonomy, we can show
that $card(F_{1}(x)\cap F_{2}(x))\leq card(F _{1}(z)\cap
F_{2}(z))$, but we have already proven that $card(F_{1}(z)\cap
F_{2}(z))\leq card(F_{1}(x_{0})\cap F_{2}(x_{0}))$.
\end{proof}

Take $x_{0}=p(a_{0},b_{0})\in M$ such that $F_{1}(x_{0})$ has not holonomy and let
$\Phi:F_{1}(x_{0})\times_{(\lambda_{1},\rho_{2})}M_{2}\rightarrow M$ be the semi-Riemannian covering map constructed in Theorem \ref{teoprin},
which has $\Omega=\Gamma/\left( \Gamma_{1}^{(a_{0},b_{0})}\times \{id\}\right)  $ as deck transformation group. Take $x\in F_{2}(x_{0})$ and a
point $b\in M_{2}$ with $\Phi
(x_{0},b)=x$. Applying Lemma \ref{coveringinducidos}, $\Phi_{1}^{(x_{0}%
,b)}:F_{1}(x_{0})\times\{b\}\rightarrow F_{1}(x)$, the restriction
of $\Phi$, is a normal semi-Riemannian covering map. Call
$\Omega_{1}^{(x_{0},b)}$ its deck transformations group.

\begin{theorem}\label{Covering hojas} In the above situation, the following sequence is exact%
\begin{equation*}
0\longrightarrow\pi_{1}(F_{1}(x_{0}),x_{0})\overset{\Phi_{1\#}^{(x_{0},b)}%
}{\longrightarrow}\pi_{1}(F_{1}(x),x)\overset{H}{\longrightarrow}%
Hol(F_{1}(x))\longrightarrow0,
\end{equation*}
where $H:\pi_{1}(F_{1}(x),x)\longrightarrow Hol(F_{1}(x))$ is the
usual holonomy homomorphism. In particular we have
$\Omega_{1}^{(x_{0},b)}=Hol\left(  F_{1}(x)\right)  $.
\end{theorem}
\begin{proof} It is clear that $\Phi_{1\#}$ is injective and $H$
is onto, so we only prove that $Ker\ H=Im\ \Phi_{1\#}$.

Take $[\gamma]\in\pi_{1}(F_{1}(x),x)$ such that $H([\gamma])=1$, i.e., $f_{\gamma}=id$, where $f_{\gamma}$ is the associated holonomy map. Take
$\alpha$ a lift of $\gamma$ in $F_{1}(x_{0})\times\{b\}$ with basepoint $(x_{0},b)$ and $\phi\in\Omega_{1}^{(x_{0},b)}$ such that $\phi
(x_{0},b)=\alpha(1)$. Since $f_\gamma=id$ it follows that $\rho_{2}\circ\phi=\rho_{2}$ and
 $\phi\times id\in\Omega$, see Remark \ref{subgrupo2}.

Therefore, taking into account that $\Phi(x,b_{0})=x$ for all $x\in F_1(x_0)$, we get
$x_{0}=\Phi_1(x_{0},b_{0})=\Phi_1(\phi(x_{0}),b_{0})=\phi(x_{0})$. Hence $\alpha$ is a loop at $x_{0}$ which holds
$\Phi_{1\#}([\alpha])=[\gamma]$. This shows that $Ker\ H\subset Im\ \Phi_{1\#}$. The other inclusion is trivial because the holonomy of the
first canonical foliation in the product $F_{1} (x_{0})\times M_{2}$ is trivial and $\Phi$ preserves the foliations.
\end{proof}

Summarizing, we obtain

\begin{corollary}\label{rechojas}
Let $M=\left(  M_{1}\times_{(\lambda_{1},\lambda_{2})}%
M_{2}\right)  /\Gamma$ be a quotient of a doubly warped product and take $x_0\in M$ such that $F_1(x_0)$ has not holonomy.
\begin{enumerate}
\item For any leaf $F_1$ there exists a normal semi-Riemannian covering map $\Phi:F _{1}(x_{0})\rightarrow F_{1}$ with deck transformation group
$Hol(F_{1})$. \item All leaves without holonomy are homothetic.
\end{enumerate}
\end{corollary}
\begin{proof}
For the first point, note that given any leaf $F_1$ it always exists $x\in F_2(x_0)$ such that $F_1=F_1(x)$.
For the second statement just note that $F_{1}(x_{0})$ and $F_{1}(x_{0}%
)\times\{b\}$ are homothetic.
\end{proof}

\begin{corollary} Let $M=\left(  M_{1}\times_{(\lambda_{1},\lambda_{2})}M_{2}\right)
/\Gamma$ be a quotient of a doubly warped product. If there is a
noncompact leaf, then any compact leaf has nontrivial holonomy.
\end{corollary}

\begin{corollary}\label{grupofundamentalgruponormal} Let $M=\left( M_{1}\times_{(\lambda _{1},\lambda_{2})}M_{2}\right)  /\Gamma$ be a quotient
of a doubly warped product and take $x_{0}=p(a_{0},b_{0})\in M$. If $F _{1}(x_{0})$ has not holonomy, then $\pi_{1}(F_{1}(x_{0}),x_{0})$ is a
normal subgroup of $\pi_{1}(M,x_{0})$.
\end{corollary}

\begin{proof} From Theorem \ref{teoprin} it is immediate that
$\pi_1(F_1(x_0),x_0)$ is a subgroup of $\pi_{1}(M,x_{0})$.
Take $[\alpha]\in\pi_{1}(F_{1}(x_{0}),x_{0})$ and $[\gamma]\in\pi_{1}%
(M,x_{0})$. We show that $[\gamma\cdot\alpha\cdot\gamma^{-1}]$ is homotopic to a loop in $F_{1}(x_{0})$. Take the covering map
$\Phi:F_{1}(x_{0})\times
M_{2}\rightarrow M$ and $\tilde{\gamma}%
=(\tilde{\gamma}_{1},\tilde{\gamma}_{2})$ a lift of $\gamma$ starting at $(x_{0},b_{0})$. Since $\Phi(\tilde{\gamma}(1))=x_{0}$, using the above
corollary we have $\Phi:F_{1}(x_{0})\times\{\tilde{\gamma}_{2}(1)\}\rightarrow F_{1}(x_{0})$ is an isometry, thus we can lift the loop $\alpha$
to a loop $\tilde{\alpha}$ starting at $\tilde{\gamma}(1)$. Therefore, the lift of $\gamma\cdot\alpha\cdot\gamma^{-1}$ to $F_{1}(x_{0})\times
M_{2}$ starting at $(x_{0},b_{0})$ is $\tilde{\gamma}\cdot\widetilde{{\alpha}}\cdot \tilde{\gamma}^{-1}$. But it is clear that this last loop is
homotopic to a loop in $F_{1}(x_{0})\times\{b_{0}\}$ and so $\gamma\cdot\alpha\cdot\gamma^{-1}$ is homotopic to a loop in $F_{1}(x_{0})$.
\end{proof}

\begin{example}\label{ejemplo} Lemma \ref{subgrupo}, and thus the results that
depend on it, does not hold if we consider more general products
than doubly warped product, as the following example shows.

First, we are going to construct a function $\lambda:\mathbb{R}^{2}%
\rightarrow\mathbb{R}^{+}$ step by step. Take $h:\mathbb{R}\rightarrow \mathbb{R}$, such that $h=id$ in a neighborhood of $0$ and
$h^{\prime}(x)>0$ for all $x\in\mathbb{R}$, and $\lambda:(-\varepsilon,\varepsilon )\times\mathbb{R}\rightarrow\mathbb{R}^{+}$ any $C^{\infty}$
function for $\varepsilon<\frac{1}{2}$. We extend $\lambda$ to the trip $(1-\varepsilon,1+\varepsilon )\times\mathbb{R}$ defining
$\lambda(x,y)=\lambda(x-1,h(y))h^{\prime}(y)$ for every $(x,y)\in (-\varepsilon,\varepsilon )\times\mathbb{R}$. Now extend it again to
$[\varepsilon,1-\varepsilon]\times\mathbb{R}$ in any way such that $\lambda:(-\varepsilon,1+\varepsilon)\times\mathbb{R}\rightarrow
\mathbb{R}^{+}$ is $C^{\infty}$. Thus, we have a function with $\lambda (x,y)=\lambda(x-1,h(y))h^{\prime}(y)$ for all $(x,y)\in(1-\varepsilon
,1+\varepsilon)\times\mathbb{R}$, or equivalently, $\lambda(x,y)=\lambda (x+1,f(y))f^{\prime}(y)$ for all $(x,y)\in(-\varepsilon,\varepsilon
)\times\mathbb{R}$, where $f$ is the inverse of $h$.

Now, we define $\lambda$ in $[1+\varepsilon,\infty)$ recursively
by $\lambda(x,y)=\lambda(x-1,h(y))h^{\prime}(y)$ and in
$(-\infty,-\varepsilon]$ by
$\lambda(x,y)=\lambda(x+1,f(y))f^{\prime}(y)$. It is easy to show
that $\lambda:\mathbb{R}^{2}\rightarrow\mathbb{R}^{+}$ is
$C^{\infty}$.

Take $\mathbb{R}^{2}$ endowed with the twisted metric
$dx^{2}+\lambda (x,y)^{2}dy^{2}$ and $\Gamma$ the group generated
by the isometry $\phi(x,y)=(x+1,f(y))$, which preserves the
canonical foliations and acts in a
properly discontinuous manner. Take $p:\mathbb{R}^{2}\rightarrow\mathbb{R}%
^{2}/\Gamma=M$ the projection. The leaf of the first foliation
through $p(0,0)$ is diffeomorphic to $\mathbb{S}^{1}$ and have not
holonomy. But
Theorem \ref{teoprin} does not hold because if $\Phi:\mathbb{S}^{1}%
\times\mathbb{R}\rightarrow M$ were a covering map, then
$\mathbb{S}^{1}$
would be a covering of all leaves of the first foliation (Corollary \ref{rechojas}%
). But this is impossible because for a suitable choice of $h$,
there are leaves diffeomorphic to $\mathbb{R}$.
\end{example}

We finish this section with a cohomological obstruction to the existence of a quotient of a doubly warped product with compact leaves. If
$M_{1}$ and $M_{2}$ are $n$-dimensional, compact and oriented manifold, K\"{u}nneth formula implies that the $n$-th Betti number of the product
$M_{1}\times M_{2}$ is greater or equal than $2$. The following theorem shows that the same is true for any oriented quotient of a doubly warped
product with $n$-dimensional compact leaves.

\begin{theorem} Let $M=\left(  M_{1}\times_{(\lambda_{1},\lambda_{2})}M_{2}\right)
/\Gamma$ be an oriented quotient of a doubly warped product such that the leaves of both foliations on $M$ are $n$-dimensional and compact
submanifolds of $M$. Then, the $n$-th Betti number of $M$ satisfies $b_{n}\geq2$.
\end{theorem}

\begin{proof} Take a point $(a_{0},b_{0})\in M_{1}\times M_{2}$ such that $x_{0}%
=p(a_{0},b_{0})$ has not holonomy and $\Phi:F_{1}(x_{0})\times_{(\rho_{1}%
,\rho_{2})}F_{2}(x_{0})\longrightarrow M$ the covering map given in Corollary \ref{corprin}. Since $M$ is oriented, $M_1$ and $M_2$ are
orientable and $\Gamma$ preserves the orientation. But $\Gamma_i^{(a_0,b_0)}$ is a normal subgroup of $\Gamma$ and therefore it preserves the
orientation of $M_i$. Thus $F_i(x_0)=M_i/\Gamma_i^{(a_0,b_0)}$ is orientable.

Let $[\varpi_{1}],[\varpi_{2}]\in H^{n}(M)$ be the Poincar\'{e} dual of $F_{1}(x_{0})$ and $F_{2}(x_{0})$ respectively. The submanifolds
$S_{i}=\Phi^{-1}(F_{i}(x_{0}))$ are closed in $F_{1}(x_{0})\times F_{2}%
(x_{0})$ and therefore they are compact. With the appropriate orientation, they have Poincar\'{e} duals
$[\sigma_{i}]=\Phi^{\ast}([\varpi_{i}])$, \cite{BottTu1982}.

Call $\pi_{i}:F_{1}(x_{0})\times F_{2}(x_{0})\rightarrow F_{i}(x_{0})$ the canonical projection, $\Phi_{i}:S_{i}\rightarrow F_{i}(x_{0})$ the
restriction of $\Phi$ to $S_{i}$, and $i_j:S_j\rightarrow F_{1}(x_{0})\times F_{2}(x_{0})$ the canonical inclusion. Consider the following
commutative diagram
\begin{align*}
\xymatrix{ S_{j} \ar[d]_{i_j} \ar[r]^{\Phi_j} & F_j(x_0) \\ F_{1}(x_{0})\times F_{2}(x_{0}) \ar[ru]^{\pi_j}}
\end{align*}

If $\Theta_{1}$ is a volume form of $F_{1}(x_{0})$, then $\Phi_{1}^{\ast }(\Theta_{1})=i_1^{\ast}(\pi_{1}^{\ast}(\Theta_{1}))$ is a volume form
in $S_{1}$. Therefore
\[
0\neq\int_{S_{1}}i_1^\ast \pi_{1}^{\ast}(\Theta_{1})=\int_{F_{1}\times F_{2}}\pi _{1}^{\ast}(\Theta_{1})\wedge\sigma_{1},
\]
thus $[\sigma_{1}]$ is not null. In the same way we can show that $[\sigma_{2}]$ is not null.

Now if $\sigma_{1}-c\sigma_{2}=d\tau$ for some $0\neq c\in\mathbb{R}$ and $\tau\in\Lambda^{n-1}(M)$, then
\[
\int_{F_{1}\times F_{2}}\pi_{1}^{\ast}(\Theta_{1})\wedge\sigma_{1}=c\int_{F_{1}\times
F_{2}}\pi_{1}^{\ast}(\Theta_{1})\wedge\sigma_{2}=c\int_{S_{2}}i_2^\ast \pi_{1}^{\ast}(\Theta_{1})=0,
\]
which is a contradiction. Therefore $[\sigma_{1}]$ and $[\sigma_{2}]$ are linearly independent, so the same is true for $[\varpi_{1}]$ and
$[\varpi_{2}]$.
\end{proof}

Observe that if the dimension of the foliations are $n$ and $m$
with $n\neq m $, then we can only conclude that the $n$-th and
$m$-th Betti numbers of $M $ satisfy $b_{n},b_{m}\geq1$. In the
category of four dimensional Lorentzian manifolds we have the
following result.

\begin{corollary} In the conditions of the above theorem, if $M$ is a four
dimensional Lorentzian manifold, then its first and second Betti
numbers satisfies $b_{1},b_{2}\geq2$.
\end{corollary}
\begin{proof} It is clear that $M$ is compact (Corollary \ref{corprin}). The
existence of a Lorentz metric implies that the Euler
characteristic is null, thus $2b_{1}=2+b_{2}$.
\end{proof}

\section{Space of leaves}

Given a foliation $\mathcal{F}$ on a manifold $M$, a point $x$ is called regular if it exists an adapted chart $(U,\varphi)$ to $\mathcal{F}$,
with $x\in U$, such that each leaf of the foliation intersects $U$ in an unique slice. The open set $U$ is also called a regular neighborhood of
$x$. If all points are regular (i.e. $\mathcal{F}$ is a regular foliation), then the space of leaves $\mathfrak{L}$ of $\mathcal{F}$ is a
manifold except for the Hausdorffness, and the canonical projection $\eta:M\rightarrow\mathfrak{L}$ is an open map, \cite{KobNom,Palais}.

Given $M=\left(M_{1}\times_{(\lambda_{1},\lambda_{2})}M_{2}\right)
/\Gamma$ a quotient of a doubly warped product, we call
$\mathfrak{L}_{i}$ the space of leaves of the induced foliations
$\mathcal{F}_{i}$ on $M$. Take $x_{0}\in M$ without holonomy and
the normal covering map $\Phi:F_{1}(x_{0})\times
_{(\rho_{1},\rho_{2})}F_{2}(x_{0})\rightarrow M$, whose group of
deck transformation is $\Psi=\Gamma/(\Gamma_{1}\times\Gamma
_{2})$. The set $\Sigma_{x_{0}}$ formed by those maps $\psi\in Diff(F_{2}%
(x_{0}))$ such that there exists $\phi\in Diff(F_{1}(x_{0}))$ with
$\phi \times\psi\in\Psi$ is a group of homotheties of
$F_{2}(x_{0})$.

\begin{lemma}\label{AccionLibre} Let $M=\left(  M_{1}\times_{(\lambda_{1},\lambda_{2})}%
M_{2}\right)  /\Gamma$ be a quotient of a doubly warped product.
Suppose that the foliation $\mathcal{F}_{1}$ has not holonomy and
take $x_{0}\in M$ such that $F_{2}(x_{0})$ has not holonomy. Then
the action of $\Sigma_{x_{0}}$ on $F_{2}(x_{0})$ is free.
\end{lemma}
\begin{proof} Take $\psi\in\Sigma_{x_{0}}$ and suppose that it has a fixed point
$x\in F_{2}(x_{0})$. If $\phi\in Diff(F_{1}(x_{0}))$ with $\phi\times\psi\in\Psi$, then $\Phi(z,x)=\Phi(\phi(z),\psi(x))=\Phi(\phi(z),x)$ for
all $z\in F_{1}(x_{0})$, but since $\mathcal{F}_{1}$ has not holonomy, applying Corollary \ref{rechojas},
$\Phi:F_{1}(x_{0})\times\{x\}\rightarrow F_{1}(x)$ is an isometry. Therefore $\phi=id$ and hence $\psi=id$.
\end{proof}

\begin{theorem}\label{Recubridor sobre hojas} Let $M=\left(  M_{1}\times_{(\lambda_{1}%
,\lambda_{2})}M_{2}\right)  /\Gamma$ be a quotient of a doubly
warped product such that $\mathcal{F}_{1}$ is a regular foliation.
If $F_{2}(x_{0})$ has not holonomy then
\begin{enumerate}
\item The group $\Sigma_{x_{0}}$ acts in a properly discontinuous
manner (in the topological sense) on $F_{2}(x_{0})$.
\item The restriction $\eta_{x_{0}}=\eta|_{F_{2}(x_{0})}:F_{2}(x_{0}%
)\rightarrow\mathfrak{L}_{1}$ is a normal covering map with
$\Sigma_{x_{0}}$ as deck transformation group.
\end{enumerate}
\end{theorem}
\begin{proof}
\begin{enumerate}
\item Suppose that $\Sigma_{x_{0}}$ does not act in a properly discontinuous manner. Then, there exists $x\in F_{2}(x_{0})$ such that for all
neighborhood $U$ of $x$ in $F_{2}(x_{0})$ there is $\psi\in\Sigma_{x_{0}}$, $\psi\neq id$, with $U\cap\psi(U)\neq\emptyset$.

Take $V\subset M$ a regular neighborhood of $x$ adapted to
$\mathcal{F}_{1}$. Since $\Phi(x_{0},x)=x$, we can lift $V$
through the covering $\Phi:F_1(x_0)\times F_2(x_0)\rightarrow M$
and suppose that there are $U_{i}\subset F_{i}(x_{0})$ open sets
with $x_{0}\in U_{1}$, $x\in U_{2}$ and $\Phi:U_{1}\times
U_{2}\rightarrow V$ an isometry. Using that $\Sigma_{x_{0}}$ does
not act in a properly discontinuous manner, there is
$\psi\in\Sigma_{x_{0}}$, $\psi\neq id$, with $y=\psi(z)$ for
certain $y,z\in U_{2}$. Moreover, $z\neq y$ since $\psi$ does not
have fixed points (Lemma \ref{AccionLibre}).

If we take $\phi$ with $\phi\times\psi\in\Psi$, then $z=\Phi(x_{0},z)=\Phi (\phi(x_{0}),y)$ and thus $F_{1}(z)=F_{1}(y)$. Now, $\Phi(U_{1}\times
\{y\})$ and $\Phi(U_{1}\times \{z\})$ are two different slices of $\mathcal{F}_{1}$ in $V$ which belong to the same leaf $F_{1}(z)$.
Contradiction.

\item It is easy to show that
$F_{2}(x_{0})/\Sigma_{x_{0}}=\mathfrak{L}_{1}$, where the
identification is $[x]\longleftrightarrow F_{1}(x)$.
\end{enumerate}
\end{proof}

\begin{corollary} Let $M=\left(  M_{1}\times_{(\lambda_{1},\lambda_{2})}M_{2}\right)
/\Gamma$ be a quotient of a doubly warped product such that
$\mathcal{F}_{1}$ is a regular foliation. If the space of leaves
$\mathfrak{L}_{1}$ is simply connected, then $M$ is isometric to a
global doubly warped product
$F_{1}\times_{(\rho_{1},\rho_{2})}F_{2}$.
\end{corollary}

We give some conditions for $\mathfrak{L}_{1}$ to be a true manifold. Recall that given $(a,b)\in F_1(x_0)\times F_2(x_0)$ with $\Phi(a,b)=x$ we
denote $\Psi_1^{(a,b)}$ the deck transformation group of the restriction $\Phi_1^{(a,b)}:F_1(x_0)\times \{b\}\rightarrow F_1(x)$ of  the
covering map $\Phi:F_1(x_0)\times F_2(x_0)\rightarrow M$.

\begin{theorem} Let $M=\left(  M_{1}\times_{(\lambda_{1},\lambda_{2})}M_{2}\right)
/\Gamma$ be a quotient of a doubly warped product. If $M_{2}$ is a complete Riemann manifold and $\mathcal{F}_{1}$ a regular foliation, then the
space of leaves $\mathfrak{L}_{1}$ is a Riemannian manifold.\end{theorem}
\begin{proof} Take $x_{0}\in M$ without holonomy. We show that $\Sigma_{x_{0}}$
is a group of isometries. Take $\psi\in\Sigma_{x_{0}}$ and $\phi:F_{1}(x_{0})\rightarrow F_{1}(x_{0})$ such that $f=\phi\times\psi\in\Psi$. As
we already said, there exist a constant $c$ such that $\psi^{\ast}(g_{2})=c^2g_{2}$ and $\rho_{2}=c\rho_{2}\circ\phi$.

Suppose $c\neq1$. Taking the inverse of $\psi$ if it were necessary, we can suppose $c<1$. Then $\psi:F_{2}(x_{0})\rightarrow F_{2}(x_{0})$ is a
contractive map and it is assured the existence of a fixed point $b\in
F_2(x_{0})$. Therefore $f(F_{1}(x_{0})\times\{b\})=F_{1}(x_{0}%
)\times\{b\}$ and $f|_{F_{1}(x_{0})\times\{b\}}\in\Psi_{1}%
^{(a,b)}$, where $a\in F_{1}(x_{0})$ is some point. Using Lemma \ref{subgrupo}, $c=1$ and we get a contradiction.

Using the above theorem, $\Sigma_{x_0}$ acts in a properly and
discontinuously
 manner on $F_2(x_0)$ in the topological sense, but since $F_2(x_0)$ is Riemannian
  and $\Sigma_{x_0}$ a group of isometries, it actually acts in a properly and discontinuously
   manner in the differentiable sense, i.e., points in different orbits have open neighborhood with disjoint orbits.
Thus, $\mathfrak{L}_{1}$ is a Riemannian manifold.
\end{proof}

\begin{remark}Given a nondegenerate foliation $\mathcal F$, it
is called semi-Riemannian (or metric) when, locally, the leaves coincide with the fibers of a semi-Riemannian submersion,
\cite{Reinhart,Walschap}. If the orthogonal distribution is integrable, then $\mathcal F$ is a semi-Riemannian foliation if and only if
$\mathcal F^{\perp}$ is totally geodesic, \cite{Whitt}. In the case of a doubly warped product $F_1(x_0)\times_{(\rho_1,\rho_2)}F_2(x_0)$, the
first canonical foliation is semi-Riemannian for the conformal metric $\bigl(\frac{\rho_1}{\rho_{2}}\bigr)^{2}g_1+g_2$.

In the hypotheses of the above theorem, $\rho_{2}$ is invariant under $\Psi$ and thus there exists a function $\sigma
_{2}:M\rightarrow\mathbb{R}^{+}$ such that $\sigma_{2}\circ\Phi=\rho_{2}$. In this case, it is easy to show that $\mathcal{F}_{1}$ is
semi-Riemannian for the conformal metric $\frac{1}{\sigma _{2}^{2}}g$, where $g$ is the induced metric on $M$. Observe that in the Riemannian
case, under regularity hypothesis, it is known that the space of leave of a Riemannian foliation is a true manifold and, moreover, the manifold
is a fiber bundle over it \cite{Reinhart}, but there is not an analogous in the semi-Riemannian case.

\end{remark}

\begin{corollary} Let $M=\left(  M_{1}\times_{(\lambda_{1},1)}M_{2}\right)
/\Gamma$ be a quotient of a warped product, where $M_2$ is a complete Riemannian manifold. If $\mathcal{F}_1$ is a regular foliation, then the
projection $\eta:M\rightarrow\mathfrak{L}_1$ is a semi-Riemannian submersion.
\end{corollary}
\begin{proof}
We already know that $\mathfrak{L}_1$ is a Riemannian manifold and $\eta_{x_0}:F_2(x_0)\rightarrow \mathfrak{L}_1$ a local isometry, where $x_0$
has not holonomy. Given $x\in F_{1}(x_{0})$, the following diagram is commutative
\begin{equation*}
\xymatrix{ \{x\}\times F_2(x_0)\ar[d]_{\eta_{x_0}\circ pr_2} \ar[r]^{\Phi} & F_2(x) \ar[ld]^{\eta_x} \\ \mathfrak{L}_1}
\end{equation*}
Since $\lambda_2=1$, the map $\Phi:\{x\}\times
F_{2}(x_{0})\rightarrow F_{2}(x)$ is a local isometry for all
$x\in F_{1}(x_{0})$. Thus, $\eta_{x}:F_{2}(x)\rightarrow
\mathfrak{L}_{1}$ is a local isometry for all $x\in M$ and
therefore, $\eta:M\rightarrow\mathfrak{L}_{1}$ is a
semi-Riemannian submersion.
\end{proof}

Observe that in the corollary, the fibres are the leaves of a warped structure, thus they are automatically umbilic.

\begin{theorem}\label{fibrebundle} Let $M=\left(  M_{1}\times_{(\lambda_{1},\lambda_{2}%
)}M_{2}\right)  /\Gamma$ be a quotient of a doubly warped product
such that $\mathcal{F}_{1}$ is a regular foliation. Then
\begin{enumerate}
\item The projection $\eta:M\rightarrow\mathfrak{L}_{1}$ is a fiber bundle. Moreover, we have $\pi_{1}(\mathfrak{L}_{1},F_{1})=\pi_{1}(M,x)/\pi
_{1}(F_{1},x)$ where $x\in F_{1}\in\mathfrak{L}_{1}$. \item There exists an open dense subset $W\subset M$ globally isometric to a doubly warped
product.
\end{enumerate}
\end{theorem}
\begin{proof}\begin{enumerate}\item Take $F_{1}\in\mathfrak{L}_{1}$ and $x_{0}\in F_{1}$ a point
without holonomy. Since
$\eta_{x_{0}}:F_{2}(x_{0})\rightarrow\mathfrak{L}_{1}$ is a
covering map, there are open sets $U\subset F_{2}(x_{0})$ and
$V\subset\mathfrak{L}_{1}$ with $x_{0}\in U$ and $F_{1}\in V$ such
that $\eta_{x_{0}}:U\rightarrow V$ is a diffeomorphism.

Now we show that $\Phi(F_{1}(x_{0})\times U)=\eta^{-1}(V)$. If
$(a,b)\in
F_{1}(x_{0})\times U$, then $\eta(\Phi(a,b))=\eta(\Phi(x_{0},b))=\eta_{x_{0}%
}(b)\in V$. Given $x\in\eta^{-1}(V)$, if we call $b=\eta_{x_{0}}%
^{-1}(\eta(x))\in U$, then $\eta(x)=\eta_{x_{0}}(b)$ and $\Phi:F_{1}%
(x_{0})\times\{b\}\rightarrow F_{1}(x)$ is an isometry because
$F_{1} $ has
not holonomy (Corollary \ref{rechojas}). Thus, there exists $a\in F_{1}%
(x_{0})$ with $\Phi(a,b)=x$.

The map $\Phi:F_{1}(x_{0})\times U\rightarrow\eta^{-1}(V)$ is
injective (and therefore a diffeomorphism). In fact, if
$(a,b),(a^{\prime},b^{\prime})\in F_{1}(x_{0})\times U$ with
$\Phi(a,b)=\Phi(a^{\prime},b^{\prime})$ then
\begin{align*}
\eta_{x_{0}}(b)  &  =\eta_{x_{0}}(\Phi(x_{0},b))=\eta_{x_{0}}(\Phi(a,b))\\
&  =\eta_{x_{0}}(\Phi(a^{\prime},b^{\prime}))=\eta_{x_{0}}(\Phi(x_{0}%
,b^{\prime}))=\eta_{x_{0}}(b^{\prime}).
\end{align*}
But since $b,b^{\prime}\in U$, we get that $b=b^{\prime}$. Now,
using that $\Phi:F_{1}(x_{0})\times\{b\}\rightarrow F_{1}(b)$ is
an isometry, we deduce that $a=a^{\prime}$.

The map $h_{V}$ that makes commutative the following diagram%
\begin{equation*}
\xymatrix{ F_1(x_0)\times U \ar[d]_{\Phi} \ar[rd]^{id\times
\eta_{x_0}} \\ \eta^{-1}(V) \ar[r]^{h_V} & F_1(x_0)\times V}
\end{equation*}
shows that $M$ is locally trivial.

Finally, using Theorem 4.41 of \cite{Hatcher}, $\eta_{\#}:\pi_{1}%
(M,F_{1},x_{0})\rightarrow\pi_{1}(\mathfrak{L}_{1},F_{1})$ is an isomorphism. But $\pi_{1}(F_{1},x_{0})$ is a normal subgroup of
$\pi_{1}(M,x_{0})$
(Corollary \ref{grupofundamentalgruponormal}), hence $\pi_{1}(M,F_{1}%
,x_{0})=\pi_{1}(M,x_{0})/\pi_{1}(F_{1},x_{0})$. \item Since $\eta_{x_0}:F_2(x_0)\rightarrow \mathfrak{L}_1$ is a covering map, we can take an
open dense set $\Theta\subset\mathfrak{L}_{1}$ and an open set $U\subset F_{2}(x_{0})$ such that $\eta_{x_{0}}:U\rightarrow\Theta$
is a diffeomorphism. Given $F_{1}\in\Theta$ we have $\Phi(x_{0},\eta_{x_{0}%
}^{-1}(F_{1}))=\eta_{x_{0}}^{-1}(F_{1})\in F_{1}$ and thus the
restriction
$\Phi:F_{1}(x_{0})\times\{\eta_{x_{0}}^{-1}(F_{1})\}\rightarrow
F_{1}$ is an isometry. Now, since $\Theta$ is dense,
$W=\eta^{-1}(\Theta)$ is dense in $M$, and taking $V=\Theta$ and
$W=\eta^{-1}(\Theta)$ in the above proof we get the result.
\end{enumerate}
\end{proof}

Recall that this open set $W$ is obtained removing a suitable set
of leaves of $\mathcal{F}_{1}$ from $M$. This is false for more
general product, as twisted products (see example \ref{ejemplo}).

\begin{remark} In \cite{Romero} a notion of local warped product on a manifold is
given as follows. Take a fiber bundle $\Pi:M\rightarrow B$ with
fibre $F$, where $M$, $B$ and $F$ are semi-Riemannian manifolds.
Suppose that there is a function
$\lambda:B\rightarrow\mathbb{R}^{+}$ such that we can take a
covering
$\{U_{i}:i\in I\}$ of trivializing open sets of $B$ with $(\Pi^{-1}%
(U_{i}),g)\rightarrow(U_{i}\times F,g_{B}+\lambda^{2}g_{F})$ an
isometry for all $i\in I$. Then it is said that $M$ is a local
warped product.

It follows that the orthogonal distribution to the fibre is integrable and using that $M$ is locally isometric to a warped product, it is easy
to show that these two foliations constitute a warped structure in the sense of definition \ref{dts}. But not all warped structures arise in
this way, since the foliation induced by the fibres of a fibre bundle has not holonomy (in fact, it is a regular foliation).
\end{remark}

\section{Global decomposition}

Given a product manifold $M_{1}\times M_{2}$, a plane
$\Pi=span(X,V)$, where $X\in TM_{1}$ and $V\in TM_{2}$, is called
a mixed plane. In this section, we show how the sign of sectional
curvature of this kind of planes determines the global
decomposition of a doubly warped structure.

\begin{lemma}\label{FuncionWarpedConstante} Let $M$ be a complete
semi-Riemannian manifold of index $\nu$. Take
$\lambda:M\rightarrow\mathbb{R}^{+}$ a smooth function and
$h_{\lambda}$ its hessian endomorphism. If $\nu<\dim M$ and
$g(h_{\lambda }(X),X)\leq0$ for all spacelike vector $X$ (or
$0<\nu$ and $g(h_{\lambda }(X),X)\leq0$ for all timelike vector
$X$), then $\lambda$ is constant.
\end{lemma}
\begin{proof} Suppose $\nu<\dim M$. Take $x\in M$ and $V\ni x$ a normal convex
neighborhood. Call $S(x)$ the set formed by the points $y\in V$
such that there exists a nonconstant spacelike geodesic inside
$V$ joining $x$ with $y $. It is obvious that $S(x)$ is an open
set for all $x\in M$ and it does not contain $x$. Let
$\gamma:\mathbb{R}\rightarrow M$ be a spacelike geodesic with
$\gamma(0)=x$. If we call $y(t)=\lambda(\gamma(t))$ then
$y^{\prime\prime}(t)\leq0$ and $y(t)>0$ for all $t\in\mathbb{R}$,
which implies that $\lambda$ is constant in $S(x)$. Take $x_{1}\in
S(x)$. In the same way, $\lambda$ is constant in $S(x_{1})$, which
is an open neighborhood of $x$. Since $x$ is arbitrary, $\lambda$
is constant. The case $0<\nu$ is similar taking timelike
geodesics.
\end{proof}

\begin{proposition} Let $M_{1}\times_{(\lambda_{1},\lambda_{2})}M_{2}$ be a doubly
warped product with $M_{1}$ and $M_{2}$ complete semi-Riemannian
manifold of index $\nu _{i}<\dim M_{i}$. If $K(\Pi)\geq0$ for all
spacelike mixed plane $\Pi$, then $\lambda_{1}$ and $\lambda_{2}$
are constant.
\end{proposition}
\begin{proof} First note that if $f\in C^\infty(M_1)$, then $g(h_{f}(X),X)=g_{1}(h_{f}^{1}(X),X)$ for all
$X\in\mathfrak{X}(M_{1})$, where $h_{f}$ is the hessian respect to the doubly warped metric $g$ and $h_{f}^{1}$ respect to $g_{1} $.

Suppose there exists a point $p\in M_{1}$ and a spacelike vector
$X_{p}\in T_{p}M_{1}$ such that $0\leq g(h_{2}(X),X)$. Given an
arbitrary spacelike
vector $V_{q}\in T_{q}M_{2}$, we have $0\leq K(X,V)+\frac{1}{\lambda_{2}%
}g(h_{2}(X),X)=-\frac{1}{\lambda_{1}}g(h_{1}(V),V)$ and applying
the above
lemma, $\lambda_{1}$ is constant. Therefore, $0\leq-\frac{1}{\lambda_{2}%
}g(h_{2}(X),X)$ for all spacelike vector $X$ and applying the above lemma again, $\lambda_{2}$ is constant too.

Suppose now the contrary case: for all spacelike vector $X\in TM_{1}$ we have $g(h_{2}(X),X)<0$. Then, the above lemma gives us that
$\lambda_{2}$ is constant. Thus $0\leq-\frac{1}{\lambda_{1}}g(h_{1}(V),V)$ for all spacelike vector $V$ and the above lemma ensures that
$\lambda_{1}$ is constant too.
\end{proof}

\begin{theorem}\label{TeoDG}Let $M=\left(M_{1} \times_{(\lambda_{1},\lambda_{2})}M_{2}\right)  /\Gamma$ be a
quotient of a doubly warped product, being $M_{1}$ a complete Riemannian manifold and $M_{2}$ a semi-Riemannian manifold with $0<\nu_{2}$.
Suppose that $\mathcal{F}_{2}$ has not holonomy, $K(\Pi)<0$ for all mixed nondegenerate plane $\Pi$ and $\lambda_{2}$ has some critical point.
Then $M$ is globally a doubly warped product.
\end{theorem}
\begin{proof} Suppose that there is a nonlightlike vector $V\in TM_{2}$
with $\varepsilon_{V}g(h_{1}(V),V)\leq0$. Given an arbitrary non
zero vector $X\in TM_{1}$, $span\{X,V\}$ is a nondegenerate
plane, thus
\begin{equation*}
-\frac{1}{\lambda_{2}}g(h_{2}(X),X)-\frac{\varepsilon_{V}}{\lambda_{1}}%
g(h_{1}(V),V)=K(X,V)<0
\end{equation*}
and therefore $0<g(h_{2}(X),X)$ for all $X\in TM_{1}$, $X\neq0$.

In the opposite case, $0<\varepsilon_{V}g(h_{1}(V),V)$ for all non
lightlike vector $V\in TM_2$. Applying Lemma
\ref{FuncionWarpedConstante}, we get that $\lambda_{1}$ is
constant, and therefore $g(h_{2}(X),X)=-\lambda_{2}K(X,V)>0$ for
all $X\in TM_{1}$, $X\neq0$. In any case $h_{2}$ is positive
definite and so $\lambda_{2}$ has exactly one critical point.

Take $x_{0}\in M$ without holonomy and the associated covering map
$\Phi
:F_{1}(x_{0})\times_{(\rho_{1},\rho_{2})}F_{2}(x_{0})\rightarrow
M$. Let $x_{1}\in F_{1}(x_{0})$ be the only critical point of
$\rho_{2}$. If $\phi\times\psi$ is a deck transformation of this
covering, then $\rho_{2}\circ \phi=c\rho_{2}$ for some constant
$c$, and it follows that $\phi(x_{1})\in F_{1}(x_{0})$ is a
critical point of $\rho_{2}$ too. Thus $\phi(x_1)=x_1$, but since
$\mathcal{F}_{2}$ has not holonomy, applying Lemma
\ref{AccionLibre}, we get $\phi\times\psi=id$. Thus $\Phi$ is an
isometry.
\end{proof}

Observe that in the conditions of the above theorem we can prove that $M=M_1\times_{(\rho_1,\lambda_2)}\left( M_2/\Gamma_2\right)$. In fact, let
$(a_0,b_0)\in M_1\times M_2$ such that $p(a_0,b_0)=x_0$. Since the points of the fibre $p_1^{-1}(x_1)$ are critical points of $\lambda_2$, being
$p_1:M_1\times \{b_0\} \rightarrow F_1(x_0)$ the covering map given in Lemma \ref{coveringinducidos}, and $\lambda_2$ has only one critical
point, it follows that $p_1$ is an isometry.
\begin{example} Kruskal space has warping function
with exactly one critical point. Thus, the last part of the above
proof shows that any quotient without holonomy is a global warped
product.
\end{example}
Now we apply the above results to semi-Riemannian submersions. We denote $\mathcal{H}$ and $\mathcal{V}$ the horizontal and vertical spaces and
$E^{v}$ (resp. $E^{h}$) will be the vertical (resp. horizontal) projections of a vector $E$.
\begin{lemma}\label{tensorT} Let $\pi:M\rightarrow B$ be a semi-Riemannian
submersion with umbilic fibres and $T$ and $A$ the O'Neill tensors
of $\pi$. Then for arbitrary $E,F\in\mathfrak{X}(M)$ and
$X\in\mathcal{H}$, it holds
\begin{enumerate}
\item $T(E,F)=g(E^{v},F^{v})N-g(N,F)E^{v}$, \item $(\nabla_{X}
T)(E,F)=g(F,A(X,E^{*}))N-g(N,F)A(X,E^{*})+g(E^{v},F^{v})\nabla_{X}
N-g(\nabla_{X} N,F)E^{v}$,
\end{enumerate}
where $N$ is the mean curvature vector field of the fibres and $E^{\ast}%
=E^{v}-E^{h}$.
\end{lemma}
\begin{proof} The first point is immediate. For the second, we have $(\nabla_{X}%
T)(E,F)=\nabla_{X}T(E,F)-T(\nabla_{X}E,F)-T(E,\nabla_{X}F)$. We
compute each term
\begin{align*}
\nabla_{X}T(E,F)  &  =\nabla_{X}(g(E^{v},F^{v})N-g(N,F)E^{v})\\
&  =\bigl(g(\nabla_{X}E^{v},F^{v})+g(E^{v},\nabla_{X}F^{v})\bigr)N+g(E^{v}%
,F^{v})\nabla_{X}N\\
&  -\bigl(g(\nabla_{X}N,F)+g(N,\nabla_{X}F)\bigr)E^{v}-g(N,F)\nabla_{X}E^{v}.\\
& \\
T(\nabla_{X}E,F)  &  =g((\nabla_{X}E)^{v},F^{v})N-g(N,F)(\nabla_{X}E)^{v}.\\
& \\
T(E,\nabla_{X}F)  &
=g(E^{v},(\nabla_{X}F)^{v})N-g(N,\nabla_{X}F)E^{v}.
\end{align*}

Rearranging terms and using that $\nabla_{X}E^{v}-(\nabla_{X}E)^{v}%
=A(X,E^{\ast})$, we obtain
\begin{align*}
(\nabla_{X}T)(E,F)  &
=\bigl(g(A(X,E^{*}),F^{v})+g(E^{v},A(X,F^{*}
))\bigr)N-g(N,F)A(X,E^{\ast})\\
&  +g(E^{v},F^{v})\nabla_{X}N-g(\nabla_{X}N,F)E^{v}.
\end{align*}

But
\begin{align*}
g(A(X,E^{*}),F^{v})+g(E^{v},A(X,F^{*}))  &  =-g(A(X,E^{h}),F^{v}%
)-g(E^{v},A(X,F^{h})) =\\
-g(A(X,E^{h}),F^{v})+g(A(X,E^{v}),F^{h})  & =
g(A(X,-E^{h}),F)+g(A(X,E^{v}),F)
=\\
g(A(X,E^{\ast}),F) &.
\end{align*}
And we obtain the result.
\end{proof}

We need to introduce the lightlike curvature of a degenerate plane in a Lorentzian manifold $(M,g)$, \cite{Harris1}. Fix a timelike and unitary
vector field $\xi$ and take a degenerate plane $\Pi=span(u,v)$, where $u$ is the unique lightlike vector in $\Pi$ with $g(u,\xi)=1$. We define
the lightlike sectional curvature of $\Pi$ as
\begin{equation*}
\mathcal{K}_{\xi}(\Pi)=\frac{g(R(v,u,u),v)}{g(v,v)}.
\end{equation*}

This sectional curvature depends on the choice of the unitary
timelike vector field $\xi$, but its sign does not change if we
choose another vector field. Thus, it makes sense to say positive
lightlike sectional curvature or negative lightlike sectional
curvature.

\begin{lemma} Let $(M,g)$ and $(B,h)$ be a Lorentzian and a Riemannian manifold respectively and  $\pi:M\rightarrow B$ a semi-Riemannian submersion with umbilic fibres. If $\xi\in\mathcal{V}$ is an
unitary timelike vector field and $\Pi=span(u,X)$ is a degenerate plane with $u\in\mathcal{V}$,
$X\in\mathcal{H}$, $g(u,u)=0$ and $g(u,\xi)=1$, then $\mathcal{K}_{\xi}%
(\Pi)=\frac{g(A(X,u),A(X,u))}{g(X,X)}$.
\end{lemma}
\begin{proof} Using the formulaes of \cite{ONeill2}, we have
\begin{align*}
g(X,X)\mathcal{K}_{\xi} (\Pi)  &  =g((\nabla_{X}T)(u,u),X)-g(T(u,X),T(u,X))\\
&  +g(A(X,u),A(X,u)).
\end{align*}

Since $u$ is lightlike, the first two terms are null by the above Lemma
\end{proof}

Given a warped product $M_{1}\times_{(1,\lambda_{2})}M_{2}$, the projection $\pi:M_{1}\times M_{2}\rightarrow M_{1}$ is a semi-Riemannian
submersion with umbilic fibres.  The following theorem assures the converse fact.

\begin{theorem} Let $M$ be a complete Lorentzian manifold, $B$ a Riemannian
manifold and $\pi:M\rightarrow B$ a semi-Riemannian submersion with umbilic fibres
of dimension greater than one and mean curvature vector $N$. If $K(\Pi)<0$ for all mixed
spacelike plane of $M$, and $N$ is closed with some zero, then $M$ is globally a warped product.
\end{theorem}

\begin{proof} By continuity, it follows that $M$ has nonpositive lightlike
curvature for all mixed degenerated plane and thus, applying the above Lemma, $A(X,u)=0$ for all $X\in\mathcal{H}$ and all lightlike
$u\in\mathcal{V}$. Therefore $A\equiv0$ and $\mathcal{H}$ is integrable and necessarily totally geodesic (see \cite{ONeill2}), which gives rise
to a warped structure, since  $N$ is closed. But being $M$ complete $M=(M_{1}\times _{(1,\lambda_{2})}M_{2})/\Gamma$ (see remark
\ref{remarkglobal}). Now, using the formulaes of Lemma \ref{curvatura} we can easily check that the curvature of a mixed plane $\Pi=span(X,V)$
is independent of the vertical vector $V$ and thus $K(\Pi)<0$ for all mixed nondegenerate plane. Finally, since $N$ has some zero, $\lambda_2$
has some critical point and applying Theorem \ref{TeoDG} we get the result.
\end{proof}

\section{Uniqueness of product decompositions}

In \cite{Esche}, the uniqueness of direct product decompositions of a nonnecessarily simply connected Riemannian manifold is studied, where the
uniqueness is understood in the following sense: a decomposition is unique if the corresponding foliations are uniquely determined. The authors
use a short generating set of the fundamental group in the sense of Gromov, which is based in the Riemannian distance. So, the techniques
employed can not be used directly in the semi-Riemannian case. In this section we apply the results of this paper to study the uniqueness
problem in the semi-Riemannian setting.

\begin{proposition} \label{subvariadad umbilica}Let $M=F_{1}\times\ldots\times F_{k}$ be a semi-Riemannian direct product and
$\mathcal{F}_{1},\ldots,\mathcal{F}_{k}$ the canonical foliations. Take $S$ an umbilic/geodesic submanifold of $M$ and suppose that there exists
$i\in\{1,\ldots,k\}$ such that $\mathcal{F}_{i}(x)\cap T_{x}S$ is a nondegenerate subspace with constant dimension for all $x\in S$. Then the
distributions $\mathcal{T}_{1}$ and $\mathcal{T}_{2}$ on $S$ determined by $\mathcal{T}_{1}(x)=\mathcal{F}_{i}(x)\cap T_{x}S$ and
$\mathcal{T}_{2}(x)=\mathcal{T}_{1}^{\perp}(x)\cap T_{x}S$ are integrable. Moreover, $\mathcal{T}_{1}$ is a regular and umbilic/geodesic
foliation and $\mathcal{T}_{2}$ is a geodesic one.
\end{proposition}

\begin{proof} It is clear that $\mathcal{T}_{1}$ is integrable. We show that $\mathcal{T}%
_{2}$ is integrable and geodesic in $S$.

Consider the  tensor $J$ given by $J(v_1,\ldots,v_i,\ldots,v_k)=(-v_1,\ldots,v_i,\ldots,-v_k)$, where $(v_1,\ldots,v_k)\in
TF_1\times\ldots\times TF_k$, and take $X,V,W\in \mathfrak{X}(S)$ with $X_{x}\in\mathcal{T}_{1}(x)$ and $V_{x},W_{x}\in\mathcal{T}_{2}(x)$ for
all $x\in S$. Since $\nabla J=0$, we have $0=(\nabla_{V}J)(X)=\nabla_{V}X-J(\nabla_{V}X)$, which means that
$\nabla_{V_{x}}X\in\mathcal{F}_{i}(x)$ since $\nabla_{V}X$ is invariant under $J$. Using that $S$ is umbilical and $X,V$ are orthogonal, we have
$\nabla_{V_{x}}X=\nabla_{V_{x}}^{S}X \in T_xS$. Therefore $\nabla_{V_{x}}X \in\mathcal{T}_{1}(x)$ for all $x\in S$.

Now, we have $g(\nabla_{V}^{S}W,X)=g(\nabla_{V}W,X)=-g(W,\nabla_{V}X)=0$. Thus, $\nabla_{V}^{S}W\in\mathcal{T}_{2}$ for all
$V,W\in\mathcal{T}_{2}$ which means that $\mathcal{T}_{2}$ is integrable and geodesic in $S$.

To see that $\mathcal T_1$ is umbilic, take $X,Y\in\mathfrak{X}(S)$ with $g(X,Y)=0$ and $X_{x},Y_{x}\in\mathcal{T}_{1}(x)$ for all $x\in S$. It
is easy to show that $\nabla_{X_{x}}Y\in\mathcal{F}_{i}(x)$ and since $S$ is umbilic, we have
\begin{equation*}
\nabla_{X_{x}}Y=\nabla_{X_{x}}^{S}Y\in\mathcal{F}_{i}(x)\cap T_{x}S=\mathcal{T}_{1}(x)
\end{equation*}
for all $x\in S$. Therefore the second fundamental form of the leaves of $\mathcal{T}_1$ inside $M$ satisfies $\mathbb{I}(X,Y)=0$ for every
couple of orthogonal vectors $X,Y\in \mathcal{T}_{1}$, which is equivalent to be umbilic submanifolds of $M$. The same argument with the second
fundamental form of $\mathcal{T}_1$ as a foliation of $S$ shows that $\mathcal{T}_1$ is an umbilic foliation of $S$. Observe that if $S$ is
geodesic it is clear that $\mathcal{T}_{1}$ is also geodesic.

Finally, we show that $\mathcal{T}_1$ is a regular foliation. Take the map $P:F_1\times\ldots\times F_k\rightarrow F_1\times\ldots\times
F_{i-1}\times F_{i+1}\times\ldots\times F_k$ given by $(x_1,\ldots ,x_k)\mapsto(x_1,\ldots,x_{i-1},x_{i+1},\ldots,x_k)$ and $i:T_2(p)\rightarrow
F_1\times\ldots\times F_k$ the canonical inclusion where $p\in S$ is a fixed point. The map $P\circ i$ is locally injective, since $Ker (P\circ
i)_{*x}=\mathcal{F}_i(x)\cap\mathcal{T}_2(x)=0$ for all $x\in T_2(p)$. Therefore, we can take a neighborhood $U\subset S$ of $p$ adapted to both
foliations $\mathcal{T}_1$ and $\mathcal{T}_2$ such that $(P\circ i)|_{V}$ is injective, being $V$ the slice of $T_2(p)$ in $U$ through $p$.
Since $P$ is constant through the leaves of $\mathcal{T}_1$, it follows that $U$ is a regular neighborhood of $p$.
\end{proof}

\begin{remark} Observe that if $S$ is geodesic then $\dim T_{x}S\cap\mathcal{F}_{i}(x)$ is constant for all $x\in S$.
\end{remark}

We say that a semi-Riemannian manifold is decomposable if it can be expressed globally as a direct product. In the contrary case it is
indecomposable.

\begin{lemma} Let $M=F_{1}\times\ldots\times F_{k}$ be a complete semi-Riemannian direct product and $\mathcal{F}_{1},\ldots,\mathcal{F}_{k}$  the canonical
foliations. Suppose $\mathcal{S}$ is a nondegenerate foliation
of dimension greater than one and invariant by parallel translation such that
$\mathcal{F}_{i}(p)\cap\mathcal{S}(p)=\{0\}$ for all $i\in\{1,\ldots,k\}$ and  some $p\in M$. Then the leaves of $\mathcal{S}$ are flat and
decomposable.
\end{lemma}

\begin{proof} Being all foliations invariant by parallel translation, the property
supposed at $p$ is in fact true at any other point of $M$. Take $x=(x_1,\ldots,x_{k})\in F_{1}\times\ldots\times F_{k}$ and suppose there is a
loop $\alpha_i :[0,1] \rightarrow F_{i}$ at $x_i$ and $v\in\mathcal{S}(x)$ such that $P_{\gamma }(v)\neq v$, where
$\gamma(t)=(x_1,\ldots,\alpha_i(t),\ldots,x_{k})$. If we decompose $v=\sum_{j=1}^{k}v_{j}\in\bigoplus_{j=1}^{k}\mathcal{F}_j(x)$, then
$P_{\gamma}(v)=P_{\gamma }(v_{i})+\sum_{j\neq i}^{k}v_{j}$ and so $0\neq v-P_{\gamma}(v)=v_{i}-P_{\gamma}(v_{i})\in
\mathcal{S}(x)\cap\mathcal{F}_{i}(x)$, which is a contradiction. Therefore, $P_{\gamma}(v)=v$ for all $v\in \mathcal{S}(x)$ and all loops
$\gamma$ of the form $\gamma(t)=(x_1,\ldots,\alpha_i(t),\ldots,x_k)$. Since $M=F_{1}\times\ldots\times F_{k}$ has the direct product metric,
$P_{\gamma}(v)=v$ for all $v\in\mathcal{S}(x)$ and an arbitrary loop $\gamma$ at $x$. In particular, the parallel translation along any loop of
a leaf $S$ is trivial. But this implies that it splits as a product of factors of the form $\mathbb{R}$ or $\mathbb{S}^{1}$.
\end{proof}

Given a curve $\gamma:[0,1] \rightarrow M$ we define $v_{\gamma}%
:[0,1] \rightarrow T_{\gamma(0)}M$ by $v_{\gamma}(t)=P_{\gamma,\gamma (0),\gamma(t)}^{-1}(\gamma^{\prime}(t))$, where $P$ is the parallel
translation. We will denote $\Omega^{M}_{p}(t_{1},\ldots,t_{m})$ the set of broken geodesics in $M$ which start at $p$ and with breaks at
$t_{i}$, where $0<t_{1}<\ldots<t_{m}<1$. If $\gamma\in\Omega^{M}_{p} (t_{1},\ldots,t_{m})$ then $v_{\gamma}$ is a piecewise constant function,
\begin{equation*}
v_{\gamma}(t)=\left\{
\begin{array}
[c]{c}%
v_{0}\text{ if }0\leq t\leq t_{1}\\
\ldots\\
v_{m}\text{ if }t_{m}\leq t\leq1\\
\end{array}
\right.
\end{equation*}
which we will denote by $(v_{0},\ldots,v_{m})$. On the other hand, if $M$ is complete, given $(v_{0},\ldots,v_{m})\in(T_{p} M)^{m+1}$ we can
construct a broken geodesic $\gamma\in\Omega^{M}_{p}(t_{1},\ldots,t_{m})$ with $v_{\gamma}\equiv(v_{0},\ldots,v_{m})$.

Now, suppose that a semi-Riemannian manifold $M$ splits as a direct product in two different manners, $M=F_{1}\times\ldots\times
F_k=S_{1}\times\ldots\times S_{k'}$. We call $\mathcal{F}_{1},\ldots,\mathcal{F}_{k}$ and $\mathcal{S}_{1},\ldots,\mathcal{S}_{k'}$ the
canonical foliations of each decomposition and $\pi_{i}:M\rightarrow F_{i}$, $\sigma_{i}:M\rightarrow S_{i}$ will be the canonical projections.

Observe that given a point $p\in M$, the leaf of $\mathcal{F}_{i}$ through $p$ is $F_{i}(p)=\{\pi_1(p)\}\times\ldots\times
F_{i}\times\ldots\{\pi_{k}(p)\}$. We will denote by $\Pi^p_i$ the projection $\Pi^p_i:M\rightarrow F_i(p)$ given by
$\Pi^p_i(x)=(\pi_1(p),\ldots,\pi_i(x),\ldots,\pi_k(p))$. Analogously, $\Sigma^p_i:M\rightarrow S_i(p)$ is given by
$\Sigma^p_i(x)=(\sigma_1(p),\ldots,\sigma_i(x),\ldots,\sigma_{k'}(p))$.

\begin{theorem}\label{unicidad} Let $M=F_{0}\times\ldots\times F_{k}$ be a complete semi-Riemannian
direct product with $F_0$ a maximal semi-euclidean factor and each $F_i$
indecomposable for $i>0$.
If $M=S_{0}\times\ldots\times S_{k'}$ is another decomposition
with $S_0$ a maximal semi-euclidean factor and each $S_j$ indecomposable for $j>0$ such that
$\mathcal{F}_{i}(p)\cap\mathcal{S}_{j}(p)$ is zero or a nondegenerate space for
some $p\in M$ and all $i,j$, then $k=k'$ and, after rearranging,
$\mathcal{F}_{i}=\mathcal{S}_i$ for all $i\in\{0,\ldots,k\}$.
\end{theorem}

\begin{proof} Fix $x\in M$ and suppose that $dim S_1(x)>1$ and $\mathcal{S}_1(x)\neq\mathcal{F}_{i}(x)$ for all $i\in\{0,\ldots,k\}$. Using the above lemma we have
that $\mathcal{S}_{1}(x)\cap\mathcal{F}_{i}(x)\neq 0$ for some $i\in\{0,\ldots,k\}$.
 Moreover, since $\mathcal{S}_1(x)\neq\mathcal{F}_{i}(x)$ it holds
 $\mathcal{S}_1(x)\cap\mathcal{F}_i(x)\neq\mathcal{S}_1(x)$ or $\mathcal{S}_1(x)\cap\mathcal{F}_i(x)\neq\mathcal{F}_i(x)$.
 We suppose the first one (the second case is similar).

Proposition \ref{subvariadad umbilica} ensures that $\mathcal{T}_{1}=\mathcal{F}_{i}\cap\mathcal{S}_{1}$ is a regular foliation and, since
$S_1(x)$ is a geodesic submanifold, $\mathcal{T}_{1}$ and $\mathcal{T}_{2}=\mathcal{T}_{1}^{\perp}\cap\mathcal{S}_{1}$ are two geodesic and
nondegenerate foliations in $S_{1}(x)$. We can choose $p\in S_1(x)$ such that the leaf $T_{2}(p)$ of $\mathcal{T}_2$ has not holonomy. We want
to show that $T _{1}(p)\cap T_{2}(p)=\{p\}$ and apply Corollary \ref{descglobal}. For this, fix an orthonormal basis in $T_{p}M$ and take a
definite positive metric such that this basis is orthonormal too. Call $|\cdot|$ its associated norm. Given
$\gamma\in\Omega_{p}^{M}(t_{1},\ldots,t_{m})$ with $v_{\gamma } \equiv(v_{0},\ldots,v_{m}) $ we call $|\gamma|=\sum_{j=0}^{m}|v_{j}|$.

Suppose  there is $q\in T_{1}(p)\cap T_{2}(p)$ with $p\neq q$. Then it exists a curve in $\Omega_{p}^{T_2(p)}(t_{1},\ldots,t_{m})$ joining $p$
and $q$ for certain $0<t_1<\ldots<t_m<1$ and so we can define
\begin{equation*}
r=\inf\{|\gamma|:\gamma\in\Omega_{p}^{T_{2}(p)}(t_{1},\ldots,t_{m})\text{ and }\gamma(1)=q\}.
\end{equation*}
We have that

\begin{itemize}

\item $r>0$. In fact, if $r=0$ then it exists $\gamma\in\Omega_p^{T_2(p)}(t_1,\ldots,t_m)$ with $\gamma(1)=q$ which lays in a neighborhood of
$p$ adapted to both foliations $\mathcal{T}_1$ and $\mathcal{T}_2$ and regular for $\mathcal{T}_1$. But since $T_1(p)=T_1(q)$ and $\gamma$ is a
curve is $T_2(p)$, the only possibility is $p=q$, which is a contradiction.

\item $r$ is a minimum. Take a sequence $\gamma_{n}\in\Omega_{p}^{T_{2}(p)} (t_{1},\ldots,t_{m})$ with
$v_{\gamma_{n}}\equiv(v_{0}^{n},\ldots,v_{m}^{n})$, $\gamma_{n}(1)=q$ and $|\gamma_{n}|\rightarrow r$. Then we can extract a
convergent subsequence of $(v_{0}^{n},\ldots,v_{m}^{n})$ to, say, $(v_{0}%
,\ldots,v_{m})$. Take $\gamma_0\in\Omega_p^{T_2(p)}(t_1,\ldots,t_m)$ with $v_{\gamma_{0}%
}\equiv(v_{0},\ldots,v_{m})$. Using the differentiable dependence of the solution respect to the initial conditions and the parameters of an
ordinary differential equation \cite[Appendix I]{KobNom}, it is easy to show that $\gamma_{0}(1)=lim_{n\rightarrow\infty}\gamma_{n}(1)=q$. Since
$|\gamma _{0}|=r$, the infimum is reached.
\end{itemize}

Now, take the map $\eta=\Sigma^p_{1}\circ\Pi^p_{i}:M\rightarrow S_{1}(p)$, which holds

\begin{itemize}
\item $\eta(T_{2}(p))\subset T_{2}(p)$, since $\eta$ takes geodesics into geodesics and $\eta_{\ast p}(\mathcal{T}_{2}(p))=\mathcal{T}_{2}(p)$.

\item $\eta(p)=p$ and $\eta(q)=q$.

\item $|\eta_{\ast p}(v)|\leq|v|$ and the equality holds if and only if $v\in\mathcal{T}_{1}(p)$.
\end{itemize}

Consider the broken geodesic $\alpha=\eta\circ\gamma_{0}\in\Omega^{T_{2}(p)}
_{p}(t_{1},\ldots,t_{m})$. Then, using that $\eta$ commutes with the parallel translation along any curve, we have $v_{\alpha}(t)=\eta_{*p}(v_{\gamma_{0}}%
(t))\equiv(\eta_{*p}(v_{0}),\ldots,\eta_{*p}(v_{m}))$, and so $|\alpha|<|\gamma_{0}|$. Since $\alpha(1)=q$ we get a contradiction.

Therefore $T_{1}(p)\cap T_{2}(p)=\{p\}$ and $S_{1}(p)$ can be decomposed as $T_{1}(p)\times T_{2}(p)$, which is a contradiction because $S_1$ is
indecomposable. The contradiction comes from supposing that $\mathcal{S}_1(x)\neq\mathcal{F}_{i}(x)$ for all $i\in\{0,\ldots,k\}$, thus it has
to hold that $\mathcal{S}_1(x)=\mathcal{F}_1(x)$ for example.
But  this means $\mathcal{S}_1=\mathcal{F}_1$.

Applying repeatedly the above
reasoning  we can eliminate the factors with dimension greater than one, except $S_0$,
in the decomposition $S_0\times\ldots \times S_{k^\prime}$, reducing the
problem to prove the uniqueness of the decomposition of a semi-Riemannian direct product
$S_0\times \mathbb{S}^1\times\ldots\times\mathbb{S}^1$,
where $S_0$ is semi-euclidean. But, in this product,
we can trivially change the metric to obtain a Riemannian direct product where we
can apply \cite{Esche}.
\end{proof}

Observe that the nondegeneracy hypothesis is redundant in the Riemannian case. On the contrary, in the semi-Riemannian case it is necessary as
the following example shows.

\begin{example} Take $L$ a complete and simply connected Lorentzian manifold with
a parallel lightlike vector field $U$, but such that $L$ can not be decomposed as a direct product, (for example a plane fronted wave,
\cite{Candela}). Take $M=L\times\mathbb{R}$ with the product metric and $X=U+\partial_{t}$. Then $X$ is a spacelike and parallel vector field
and since $M$ is complete and simply connected, $M$ splits as a direct product with the integral curves of $X$ as a factor. Thus $M$ admits two
different decomposition as direct product, although $L$ is indecomposable.
\end{example}

\bibliographystyle{amsplain}

\begin{thebibliography}{99}
\bibitem{Arouche} A. Arouche, M. Deffaf and A. Zeghib. On Lorentz dynamics: from group actions to warped products via homogeneous
spaces. \textit{Trans. Amer. Math. Soc.} {\bf 359} (2007), 1253-1263.

\bibitem {BottTu1982}R. Bott and L. W. Tu. \textit{Differential forms in
Algebraic Topology}. Springer, (1982).

\bibitem {Camacho}C. Camacho and A. L. Neto. \textit{Geometric theory of
foliations}, Birkh\"auser. Boston, (1985).

\bibitem {Candela}A. M Candela, J. L. Flores and M. S\'anchez. On General plane fronted waves. Geodesics. \textit{Gen. Rel. Grav.} {\bf35} (2003),
631-649.

\bibitem {Derdzinski}A. Derdzinski. Classification of certain compact
Riemannian manifolds with harmonic curvature and non-parallel Ricci tensor. \textit{Math. Z.} {\bf 172} (1980), 273-280.

\bibitem {Tischler}D. B. A. Epstein, K. C. Millett and D. Tischler.
Leaves without holonomy. \textit{J. London Math. Soc. (2)} {\bf 16} (1977), 548-552.

\bibitem {Esche}J. H. Eschenburg and E. Heintze. Unique decomposition
of Riemannian manifolds. \textit{Proc. Amer. Math. Soc.} {\bf 126} (1998), 3075-3078.

\bibitem{Escobales} R. H. Escobales Jr. and P. E. Parker. Geometric consecuences of the normal curvature cohomoly class in umbilic foliations.
 \textit{Indiana Univ. Math. J.} {\bf 37} (1988), 389-408.

\bibitem {Jelonek}W. Jelonek. Killing tensors and warped
products. \textit{Ann. Polon. Math.} {\bf 75} (2000), 15-33.

\bibitem {GutOle09}M. Guti\' errez and B. Olea. Global decomposition of a
manifold as a Generalized Robertson-Walker space. \textit{Differential Geom. Appl.} {\bf 27} (2009), 146-156.

\bibitem {Harris1}S. G. Harris. A triangle comparison theorem for
Lorentz manifolds. \textit{Indiana Univ. Math. J.} {\bf 31} (1982), 289-308.

\bibitem {Hatcher}A. Hatcher. \textit{Algebraic Topology},
http://www.math.cornell.edu/\symbol{126}hatcher/AT/ATpage.html.

\bibitem {Whitt}D. L. Johnson and L. B. Whitt. Totally geodesic
foliations \textit{J. Differential Geometry} {\bf 15} (1980), 225-235.

\bibitem {KobNom}S. Kobayashi and K. Nomizu. \textit{Foundations of
differential geometry, vol. I.} Interscience Publishers, New York,
(1963).

\bibitem {Koike}N. Koike. Totally umbilic foliations and decomposition
theorems. \textit{Saitama Math. J.} {\bf 8} (1990), 1-18.

\bibitem{Lappas} D. Lappas. Locally warped products arising from certain group actions. \textit{Japan J. Math.} {\bf 20} (1994), 365-371.

\bibitem {ONeill}B. O'Neill. \textit{Semi-Riemannian geometry with
applications to relativity}, Academic Press, New York, (1983).

\bibitem {ONeill2}B. O'Neill. The fundamental equations of a
submersion. \textit{Michigan Math. J.} {\bf 13} (1966), 459-469.

\bibitem {Palais}R. Palais, \textit{A global formulation of the Lie theory of
transformation groups}, Memoirs of the Amer. Math. Soc. {\bf 22} (1957).

\bibitem {Ponge}R. Ponge and H. Reckziegel. Twisted product in
pseudo-Riemannian geometry. \textit{Geom. Dedicata} {\bf 48} (1993), 15-25.

\bibitem {Reinhart}B. L. Reinhart. Foliated manifolds with bundle-like
metrics. \textit{Ann. of Math} {\bf 69} (1959), 119-132.

\bibitem {Romero}A. Romero and M. S\'anchez. On completeness of
certain families of semi-Riemannian manifolds. \textit{Lett. Math. Phys.} {\bf 53} (1994), 103-107.

\bibitem {Walschap}G. Walschap. Spacelike metric foliations. \textit{J. Geom.
Phys.} {\bf 32} (1999), 97-101

\bibitem {Wang}P. Wang. Decomposition theorems of Riemannian
manifold. \textit{Trans. Amer. Math. Soc.} {\bf 184} (1973), 327-341.

\bibitem {ShaZhuIgo}Y. L. Shapiro, N. I. Zhukova and V. A. Igoshin,
Fibering on some classes of Riemannian manifolds. \textit{Izv. Vyssh. Uchebn. Zaved. Mat.} {\bf 7} (1979), 93-96.
\end{thebibliography}

\end{document}